\documentclass[lettersize,journal]{IEEEtran}
\usepackage{cite}
\usepackage{amsmath,amssymb,amsfonts}
\usepackage{algorithm}
\usepackage{algcompatible}
\usepackage{graphicx}
\usepackage{graphicx}  
\usepackage{subfigure} 
\usepackage{textcomp}
\usepackage{float}
\usepackage{multicol}
\usepackage{arydshln}
\usepackage{amsmath}
\usepackage{amsthm}
\usepackage{amssymb}
\usepackage{amsfonts}
\usepackage{graphicx}
\usepackage{mathrsfs}
\usepackage{subfigure}
\usepackage{url}
\usepackage{booktabs}
\usepackage{float}
\usepackage{color}
\usepackage{bm}
\usepackage{graphicx}
\usepackage{epstopdf}
\usepackage{algpseudocode}
\usepackage{arydshln}

\newcommand{\col}{\hbox{col}}
\newtheorem{assmp}{\bf Assumption}

\newtheorem{pro}{\bf Problem}

\newtheorem{rem}{\bf Remark}

\newtheorem{lem}{\bf Lemma}

\newtheorem{thm}{\bf Theorem}

\newenvironment{Prf}{\noindent{\emph{Proof:}}}{\hfill $\Box$\par}

\newcommand{\EQ}{\begin{eqnarray}}
\newcommand{\EN}{\end{eqnarray}}
\newcommand{\EQQ}{\begin{eqnarray*}}
\newcommand{\ENN}{\end{eqnarray*}}

\ifCLASSINFOpdf
\else
\fi

\begin{document}
\title{Data-Driven Output-Based Approach to the  Output Regulation Problem of Unknown Linear Systems via Value Iteration}

\author{Haoyan~Lin~and~Jie~Huang,~\IEEEmembership{Life Fellow,~IEEE}
\thanks{This work was supported by the Research Grants Council of the Hong Kong Special Administrative Region under grant No. 14203924.}
\thanks{The authors are with the Department of Mechanical and Automation Engineering, The Chinese University of Hong Kong, Hong Kong (e-mail: hylin@mae.cuhk.edu.hk; jhuang@mae.cuhk.edu.hk. Corresponding author: Jie Huang.)}}

\maketitle

\begin{abstract}
The output regulation problem for unknown linear systems has been studied using state-based and output-based internal model approaches in the special case with no disturbances.
This paper further investigates the output regulation problem for unknown linear systems using a data-driven output-based approach via value iteration. For this purpose, we first develop a novel output-feedback control law that does not explicitly rely on the observer gain to solve the output regulation problem.
We then show that the data-driven approach for designing an output-feedback control law for the given plant can be reduced to the data-driven design of a state-feedback control law for a well-defined augmented auxiliary  system.
As a result, we develop a systematic data-driven approach to solve the output regulation problem for unknown linear systems via value iteration.
Finally, we establish a relation between the data-driven state-feedback control law and the data-driven output-feedback control law in the LQR sense.
\end{abstract}


\IEEEpeerreviewmaketitle

% has proven to be effective in addressing a wide range of optimal control problems

\section{Introduction}
As a data-driven technique, reinforcement learning (RL), also known as adaptive dynamic programming (ADP), has been developed to approximately solve various optimal control problems for systems with partially or completely unknown dynamics \cite{Bertsekas2019, Lewis2012, Meyn2022, sutton2018reinforcement, Werbos}.
Reference \cite{Vrabie2009} first studied the linear quadratic regulator (LQR) problem for linear systems with an unknown system matrix by iteratively updating the control policy, an approach known as the policy iteration (PI) method. Reference \cite{Jiang2012} further investigated the LQR problem for linear systems with both unknown system and input matrices. However, since the PI method requires an initially stabilizing feedback gain—which may be difficult to obtain when the system is unknown, reference \cite{Bian2016} proposed a value iteration (VI) method that iteratively updates the solution to a Riccati equation, thereby eliminating the need to start with an initially stabilizing feedback gain. All the references mentioned above assumed the availability of the full system state.
When only the system output is available, reference \cite{Zhu2014} considered solving the output-feedback LQR problem with a known input matrix, assuming the system is stabilizable via static output-feedback. Reference \cite{Modares2016} further addressed the output-feedback LQR problem for completely unknown linear systems by introducing a discounted cost function, where the discount factor has an upper bound determined by the system model to ensure stability. To avoid relying on a discounted cost function, references \cite{Rizvi2023} and \cite{Rizvi2019} developed both PI-based and VI-based output-feedback ADP algorithms that involve constructing a model-free observer.
The approaches in \cite{Rizvi2023} and \cite{Rizvi2019} are particularly attractive because they preserve the well-known separation principle. However, these methods require the parameterization matrix to be of full row rank. Since this matrix is unknown, verifying its rank is challenging. Recently, reference \cite{Chen2022} showed that this matrix is of full row rank if the plant is controllable, but this result does not apply when the plant is only stabilizable.
More recently, reference \cite{Lin2024} proposed a new method that successfully transforms the output-feedback LQR problem into a state-feedback LQR problem for an auxiliary  linear system derived from the model-free observer in \cite{Rizvi2019}. Compared to \cite{Rizvi2019}, the method in \cite{Lin2024} eliminates the full-row-rank condition on the parameterization matrix and reduces computational cost, as the input matrix of the auxiliary  system is known.

A more general control problem involves asymptotic tracking of a class of reference inputs and rejection of a class of disturbances. When both the references and disturbances are generated by a linear autonomous system known as the exosystem, the problem is termed the output regulation problem. This problem is more complex than mere stabilization and has been extensively studied (e.g., in \cite{Davision1976, Francis1976, Huang2004}). It can be addressed via two main approaches: feedforward control and the internal model principle. The latter converts the output regulation problem of the original system into a stabilization problem for an augmented system comprising the plant and a dynamic compensator (the internal model). This approach has the advantage of not requiring the solution of {\color{black} the regulator equations—a set of linear matrix equations.}
Recently, significant efforts have been devoted to solving the linear output regulation problem (LORP) for unknown or partially unknown systems using data-driven methods. Reference \cite{Gao2016} addressed the problem with unknown system, input, and disturbance matrices using the feedforward control approach, which requires solving the regulator equations. Subsequently, \cite{Liu2018} considered the LORP for single-input single-output systems via the internal model approach, allowing the output matrix to be unknown. Reference \cite{Xie2022} modified the approach of \cite{Liu2018} to reduce data storage requirements. More recently, \cite{LH2024} investigated the problem using the internal model approach for multi-input multi-output systems and proposed a new learning algorithm that further reduces computational cost and relaxes the solvability conditions.
The aforementioned references all assumed full state information. More recent work has focused on solving the LORP for unknown systems via output-feedback control \cite{Chen2022, Xie2023}. Specifically, \cite{Chen2022} studied a special case of the LORP with no disturbances using an output-feedback control law based on the PI method. Furthermore, \cite{Xie2023} addressed the general case using the VI method. However, the control law in \cite{Xie2023} depends on the exosystem state and is therefore not truly output-feedback. Both \cite{Chen2022} and \cite{Xie2023} employ the parameterized observer from \cite{Rizvi2023}, thus requiring the parameterization matrix to be of full row rank—a condition guaranteed only if the plant is controllable.
In this paper, we propose a novel model-free output-feedback control law for solving the LORP. The main contributions of this paper are summarized as follows:
\begin{itemize}
	\item[1)] We first develop a novel output-feedback control law that does not explicitly depend on the observer gain for solving the output regulation problem. This control law enables the design of a data-driven output-feedback method to solve the output regulation problem.
	\item[2)] We define an augmented auxiliary  system and show that solving the state-feedback stabilization problem for this system yields an output-feedback solution to the output regulation problem for the original system.
	\item[3)] Our approach does not require the parameterization matrix to be of full row rank, thereby accommodating cases where the plant is only stabilizable.
	\item[4)] By exploiting the fact that the input matrix of the augmented auxiliary  system is known, we develop a data-driven approach with significantly fewer unknown parameters. This not only substantially reduces computational cost but also relaxes the solvability conditions.
	\item[5)] We establish a relationship between the data-driven state-feedback LQR solution for the augmented auxiliary  system and the data-driven output-feedback LQR solution for the augmented system, showing that the former converges asymptotically to the latter in the sense elaborated in Theorem \ref{thm4}.
\end{itemize}

The remainder of the paper is organized as follows. Section \ref{Preliminaries} reviews both state-feedback and output-feedback approaches for solving the LQR problem of unknown linear systems via value iteration. Section \ref{regulation} presents a novel output-feedback control law that does not explicitly depend on the observer gain for solving the output regulation problem. Section \ref{main results} establishes the main results, including a systematic data-driven output-feedback control approach via value iteration for solving the output regulation problem {\color{black}of unknown linear systems.} In Section \ref{relaOP}, we further explore the relationship between the data-driven state-feedback and output-feedback control laws in the context of the LQR problem. Section \ref{eg} provides two simulation examples to demonstrate the effectiveness of the proposed method. Finally, Section \ref{conclusion} concludes the paper.

\indent\textbf{Notation}
Throughout this paper, $\mathbb{R}, \mathbb{N},$ and $\mathbb{N}_+$ represent the sets of real numbers, nonnegative integers, positive integers, respectively. %$\mathcal{S}^n$ is the set of all $n\times n$ real, symmetric matrices.
$\mathcal{P}^n$ is the set of all $n\times n$ real, symmetric and positive semidefinite matrices. $||\cdot||$ represents the Euclidean norm for vectors and the induced norm for matrices.
For $b=[b_1, b_2, \cdots, b_n]^T\in \mathbb{R}^n$, $\text{vecv}(b)=[b_1^2,b_1b_2,\cdots,b_1b_n,b_2^2,b_2b_3,\cdots, b_{n-1}b_n,b_n^2]^T \in \mathbb{R}^{\frac{n(n+1)}{2}}$. For a symmetric matrix $P=[p_{ij}]_{n\times n}\in \mathbb{R}^{n\times n}$, $\text{vecs}(P)=[p_{11},2p_{12},\cdots,2p_{1n},p_{22},2p_{23},\cdots, 2p_{n-1,n}, p_{nn}]^T\in \mathbb{R}^{\frac{n(n+1)}{2}}$. For  $v\in \mathbb{R}^n$, $|v|_P=v^TPv$. For column vectors $a_i, i=1,\cdots,s$,  $\mbox{col} (a_1,\cdots,a_s )= [a_1^T,\cdots,a_s^T  ]^T,$ and, if  $A = (a_1,\cdots,a_s )$, then  vec$(A)=\mbox{col} (a_1,\cdots,a_s )$.
For $A\in \mathbb{R}^{n\times n}$, {\color{black}$\lambda(A) $ denotes the set composed of all the eigenvalues of $A$.}  `blockdiag' denotes the block diagonal matrix operator. $I_n $ denotes the identity matrix of dimension $n$. For a positive semidefinite matrix $Q\in \mathcal{P}^n$, $\sqrt{Q}$ denotes the positive semi-definite square root.

\section{Preliminaries and Problem Formulation}\label{Preliminaries}
\indent In this section, we first review the state-based and output-based VI approach for solving the model-free LQR problem based on \cite{Bian2016} and \cite{Lin2024}, respectively.

\subsection{State-Based {\color{black}Value Iteration Approach} for Solving LQR Problem without Knowing $A$, $B$ and $C$ \cite{Bian2016}}\label{stateLQR}
Consider the following linear system:
\begin{equation}\label{lisys}
\begin{aligned}
\dot{x}=&Ax+Bu\\
y=&Cx
\end{aligned}
\end{equation}
where $x\in \mathbb{R}^n$ is the system state, $u\in\mathbb{R}^m$ is the input and $y\in \mathbb{R}^p$ is the measurement output. We make the following assumptions:
\begin{assmp}\label{ass1}
The pair $(A,B)$ is stabilizable.\end{assmp}

\begin{assmp}\label{ass2}
The pair $(A,C)$ is observable.
\end{assmp}

The  LQR problem for \eqref{lisys} is to find a control law $u=Kx$ such that the cost $\int_{0}^{\infty}(y^T{Q_y}y+u^T{R}u)d\tau$ is minimized,
where ${Q_y}={Q_y}^T\geq0, {R}={R}^T>0$, with $(A,\sqrt{{Q_y}}C)$ observable.

Let $Q=C^TQ_yC$. By \cite{kucera1972}, under Assumption \ref{ass1},  the following algebraic Riccati equation
\begin{align}\label{areli}
 A^T{P}^*+{P}^*A+{Q}-{P}^*B{R}^{-1}B^T{P}^*=0
\end{align}
admits a unique positive definite solution ${P}^*> 0$. Then, the solution to the LQR problem of \eqref{lisys} is given by $u^*=K^*x$ with ${K}^*=-{R}^{-1}B^T{P}^*$.

When $A,B,C$ are unknown, \cite{Bian2016} developed the following iterative approach to approximately solve the algebraic Riccati equation \eqref{areli}.

For convenience,  for any vectors {\color{black}$a \in \mathbb{R}^n$, $b \in \mathbb{R}^m$} and any integer $s\in \mathbb{N}_+$, define
\begin{equation}\label{defi}
\begin{aligned}
\delta _a=&[\text{vecv}(a(t_1))-\text{vecv}(a(t_0)), \cdots ,\\ &\text{vecv}(a(t_s))-\text{vecv}(a(t_{s-1}))]^T\\
I_{aa}=&[\int_{t_0}^{t_1}\text{vecv}(a)d\tau , \int_{t_1}^{t_2}\text{vecv}(a) d\tau, \cdots , \int_{t_{s-1}}^{t_s}\text{vecv}(a) d\tau]^T\\
I_{ab}=&[\int_{t_0}^{t_1}a\otimes {R}b d\tau , \int_{t_1}^{t_2}a\otimes {R}b d\tau, \cdots , \int_{t_{s-1}}^{t_s}a\otimes {R}b d\tau]^T\\
\end{aligned}
\end{equation}

{\color{black}Let ${H}_k=A^T{P}_k+{P}_kA, {K}_k=-{R}^{-1}B^T{P}_k$.
Then, integrating $\frac{ d (x^T (t) {{P}_k} x (t))}{dt}$ along the solution of  \eqref{lisys} gives
\begin{align}\label{intvipre}
|x(t+\delta t)|_{{P}_k}-|x(t)|_{{P}_k} =\int_{t}^{t+\delta t}\lbrack |x|_{{H}_k}-2{u}^T{R}{K}_kx \rbrack d\tau
\end{align}
\eqref{defi} and  \eqref{intvipre}  imply
   \begin{align}\label{vilinearpre}
   	\Psi\begin{bmatrix}
   		\text{vecs}({H}_k)\\
   		\text{vec}({K}_k)
   	\end{bmatrix}=\Phi_k
   \end{align}
   where $\Psi=\begin{bmatrix}
   	I_{xx}&-2I_{xu}
   \end{bmatrix}$ and $\Phi_k=\delta_{x}\text{vecs}({P}_k)$.}

The solvability of \eqref{vilinearpre} is guaranteed by the following lemma \cite{Bian2016}.
\begin{lem}
The matrix $\Psi$ has full column rank  if
\begin{align}\label{rankconvi}
\textup{rank}([I_{xx}, I_{xu}])=& \frac{n(n+1)}{2}+mn
\end{align}
\end{lem}

Let $\epsilon_k$ be a series of time steps satisfying
\begin{align}
   	\epsilon_k>0, \; \sum_{k=0}^{\infty}\epsilon_k=\infty, \; \sum_{k=0}^{\infty}\epsilon_k^2<\infty
\end{align}
$\{B_j\}_{j=0}^{\infty}$ be a collection of bounded subsets in $\mathcal{P}^{n}$ satisfying
\begin{align}
   	B_j\subseteq B_{j+1}, \; j\in \mathbb{N}, \; \lim\limits_{j\to \infty}B_j=\mathcal{P}^{n}
\end{align}
and let $\varepsilon>0$  be some small real number for determining the convergence criterion.
Then the solution to \eqref{areli} can be iteratively obtained  by  Algorithm \ref{vialg1}.

\begin{algorithm}
\caption{Model-free VI Algorithm for Solving \eqref{areli} \cite{Bian2016}} \label{vialg1}
\begin{algorithmic}[1]
\State Choose ${P}_0=({P}_0)^T>0$. $k,j \gets 0$.
\Loop
\State Substitute ${P}_k$ into \eqref{vilinearpre} to solve ${H}_k$ and ${K}_k$.
\State \begin{equation}\label{A2}
\begin{aligned}
\tilde{P}_{k+1} \gets {P}_k+\epsilon_k({H}_k+Q-({K}_k)^TR{K}_k)\\
\end{aligned}
\end{equation}
\If{$\tilde{P}_{k+1} \notin B_j $}
    \State ${P}_{k+1} \gets {P}_0$. $j \gets j+1$.
\ElsIf{$||\tilde{P}_{k+1}-{P}_k||/\epsilon_k<\varepsilon$}
    \State \textbf{return} ${P}_k$ and ${K}_k$ as approximations to ${P}^*$ and ${K}^*$, respectively.
\Else
	\State ${P}_{k+1} \gets \tilde{P}_{k+1}$
\EndIf
\State $k \gets k+1$
\EndLoop
\end{algorithmic}
\end{algorithm}

\subsection{Output-Based Value Iteration Approach for Solving LQR Problem without Knowing $A$, $B$ and $C$\cite{Lin2024}}\label{outputLQR}

 The approach in \cite{Bian2016} assumed that the full state of the plant is available for feedback. In most applications, the full state is not available. Reference \cite{Rizvi2019} further considered the LQR problem with only output $y$ available for feedback by developing a parameterized approach for estimating the state.
 Reference \cite{Lin2024} further modified the approach of \cite{Rizvi2019} by converting the output-based LQR problem to a state-based LQR problem of an auxiliary  system.
To summarize the approach of \cite{Rizvi2019}  and \cite{Lin2024}, we first note that, under Assumption \ref{ass2}, there exists an observer gain $L$ such that $A-LC$ is Hurwitz. Then the following Luenberger observer:
\begin{align}\label{observer}
\dot{\hat{x}}= &A\hat{x}+Bu+L(y-C\hat{x})\notag \\
=&(A-LC)\hat{x}+Bu+Ly
\end{align}
will drive the state $\hat{x}$  to $x$ exponentially.

Let $\Lambda(s)=\textup{det}(sI_n-A+LC) = s^n + \alpha_{n-1} s^{n-1} + \cdots+ \alpha_1 s + \alpha_0$. Then $(sI_n-A+LC)^{-1}=\frac{D_{n-1}s^{n-1}+D_{n-2}s^{n-2}+\cdots+D_1s+D_0}{\Lambda(s) }$ where $D_i\in \mathbb{R}^{n\times n}, i=0,1,\cdots, n-1$.
Let $\zeta^i(t)\in \mathbb{R}^n$ be governed by the following system:
\begin{align} \label{observer2}
	\dot{\zeta}^i (t)=\mathcal{A} \zeta^i (t) +b w_i(t), ~~i = 1, \cdots, m+p
\end{align}
where $w_i (t) = u_i (t)$, $i = 1, \cdots, m$,  $w_i (t) = y_{i-m}(t)$, $i = m+1, \cdots, m+p$, and
 \begin{align*}
	\mathcal{A}=\begin{bmatrix}
		0&1&0& \cdots& 0\\
		0& 0& 1& \cdots & 0\\
		\vdots & \vdots&\vdots& \vdots& \vdots\\
		0 & 0& 0& \cdots& 1\\
		-\alpha_0 & -\alpha_1& -\alpha_2& \cdots& -\alpha_{n-1}
	\end{bmatrix},~ b=\begin{bmatrix}
		0\\0\\ \vdots \\0\\ 1
	\end{bmatrix}
\end{align*}

 \begin{rem}
 Since $\Lambda(s)$ is user-defined, $\mathcal{A}$ and $b$ are known matrices.
\end{rem}

Let $\zeta = \col (\zeta^1, \cdots,  \zeta^{p+m} )$ and  $M = [M_1,\cdots ,M_{p+m}]$ where $M_i =  [D_0 f_i ,\cdots ,D_{n-1} f_i]$ with $f_i$ the $i^{th}$ column of $B$ for $i =1, \cdots, m$ and $(i - m)^{th}$ column of $L$ for $i = m+1, \cdots, m+p$. {\color{black}Since, by \eqref{observer} and \eqref{observer2},  $ \hat{\bf{x}}(s)
= (sI_n - A+ LC)^{-1} \left ([B~~ L] \begin{bmatrix}
u(s)\\y(s)
\end{bmatrix} + \hat{x} (0) \right) $ and $ M  \zeta (s)
= (sI_n - A+ LC)^{-1} [B~~ L]\begin{bmatrix}
u(s)\\y(s)
\end{bmatrix} $  where $\hat{\bf{x}}(s), \zeta (s),u(s),y(s)$ are the Laplace transform of $\hat{x}(t), \zeta (t),u(t),y(t)$, respectively, we have
 $\lim_{t \rightarrow \infty} (M\zeta (t) - \hat{x}(t)) = 0$ exponentially. Thus, $\lim_{t \rightarrow \infty} (M\zeta (t) - x(t)) = 0$ exponentially.}

By definition, the dynamics of $\zeta$ can be written in the following form:
\begin{align}\label{dyz}
	\dot{\zeta}= (I_{m+p}\otimes \mathcal{A}) \zeta+ B_\zeta u+E_\zeta y
\end{align}
where $B_\zeta = \begin{bmatrix}
	I_m\otimes b\\ {\bf0}
\end{bmatrix} $  and $E_\zeta=\begin{bmatrix}
	{\bf0}\\I_p\otimes b
\end{bmatrix}$.

{\color{black}Let $e_\zeta\triangleq M\zeta-x$. Substituting $y=Cx=CM\zeta-Ce_\zeta$ into \eqref{dyz} gives
\begin{align}\label{dyz2bf}
	\dot{\zeta}=A_\zeta\zeta+B_\zeta u- E_\zeta C e_\zeta
\end{align}
where
$A_\zeta=(I_{m+p}\otimes \mathcal{A})+\begin{bmatrix}
	{\bf0}\\I_p\otimes b
\end{bmatrix}CM$.

Since $e_\zeta$ tends to 0 exponentially, the effect of $e_\zeta$ can be ignored after some finite time $t_0>0$. As a result, \eqref{dyz2bf} can be simplified to the following form:
\begin{align}\label{dyz2}
	\dot{\zeta}=A_\zeta\zeta+B_\zeta u, ~~ t\geq t_0
\end{align}}

Note that $A_\zeta$ is an unknown matrix since $CM$ is unknown, but $B_\zeta$ is known.

It was shown in Lemmas 6 and 7 of \cite{Lin2024} that, under Assumptions \ref{ass1} and \ref{ass2}, the pair $(A_\zeta, B_\zeta)$ is stabilizable, and,  for any $Q_y>0$,  the pair $(A_\zeta, \sqrt{{Q_y}}CM)$ is detectable. Thus,
the following algebraic Riccati equation \eqref{areliz} admits a unique positive semidefinite solution ${P}_\zeta^*$:
\begin{align}\label{areliz}
	A_\zeta^T{P}_\zeta^*+{P}_\zeta^*A_\zeta+Q_\zeta-{P}_\zeta^*B_\zeta R^{-1}B_\zeta^T{P}_\zeta^*=0
\end{align}
where $Q_\zeta=M^TQM=M^TC^TQ_yCM\geq0$, $ Q_y$ and ${R}$ are the same as those of \eqref{areli}.
Further, it was shown in  Theorem 2 of \cite{Lin2024} that ${P}_\zeta^*=M^TP^*M\geq 0$. As a result, the solution to the following LQR problem:
\begin{equation}\label{op1}
		\begin{aligned}
			& \min_{u}\int_0^\infty(|\zeta|_{{Q}_\zeta}+|{u}|_R)dt, {Q}_\zeta \geq0, R>0\\
			& \text{subject to}\quad \eqref{dyz2}
		\end{aligned}
	\end{equation}
is given by $u^*_{\zeta}=K_\zeta^* \zeta$ where $K_\zeta^*=-{R}^{-1}B_\rho^T{P}_\zeta^*$. Moreover, it can be shown that $K_\zeta^*=-{R}^{-1}B^TP^*M=K^*M$ with $P^*$
 the unique positive definite solution to
\eqref{areli}. Since $M\zeta$ converges to $x$ exponentially, $u^*_\zeta=K_\zeta^*\zeta=K^*M\zeta$ converges to the optimal controller $u^*=K^*x$ for \eqref{lisys} exponentially. Therefore,  the output-feedback LQR problem of \eqref{lisys} has been converted into the state-feedback LQR problem of the auxiliary  system \eqref{dyz2} with a known input matrix.

{\color{black}Let $\bar{H}_k=A_\zeta^T\bar{P}_k+\bar{P}_kA_\zeta+{Q}_\zeta$,  $ \bar{K}_k=-{R}^{-1}B_\zeta^T\bar{P}_k$.
Then, from \eqref{dyz2}, using the relation $\zeta^TQ_\zeta\zeta=y^TQ_yy$, we obtain
\begin{align}\label{intviz}
|\zeta(t+\delta t)|_{\bar{P}_k}-&|\zeta(t)|_{\bar{P}_k} +\int_{t}^{t+\delta t} |y|_{Q_y} d\tau\notag \\=&\int_{t}^{t+\delta t}\lbrack |\zeta|_{\bar{H}_k} +2{u}^T B_\zeta^T\bar{P}_k \zeta \rbrack d\tau
\end{align}
Noting that ${R} \bar{K}_{k}=- B_\zeta^T\bar{P}_{k}$ is known, \eqref{defi} and  \eqref{intviz}  imply
\begin{align}\label{vilinearznew}
	\tilde{\Psi}\text{vecs}(\bar{H}_k)=\tilde{\Phi}_k
\end{align}
where $\tilde{\Psi}=I_{\zeta\zeta}$ and $\tilde{\Phi}_k=\delta_{\zeta}\text{vecs}(\bar{P}_k)+I_{yy}\text{vecs}(Q_y)+2I_{\zeta u}	\text{vec}(\bar{K}_k)$.}

The solvability of \eqref{vilinearznew} is guaranteed by the following condition.
\begin{lem}
	The matrix $\tilde{\Psi}$ has full column rank if
	\begin{align}\label{rankconviz}
		\textup{rank}(I_{\zeta\zeta})=& \frac{n_\zeta(n_\zeta+1)}{2}
	\end{align}
	where $n_\zeta=(m+p)n$.
\end{lem}

The output-based VI Algorithm of \cite{Lin2024} is summarized in Algorithm \ref{alg2}, where $\epsilon_k$ and $\varepsilon$ are defined in the same way as those in Section \ref{stateLQR},  $\{\bar{B}_j\}_{j=0}^{\infty}$ is a collection of bounded subsets in $\mathcal{P}^{n_\zeta}$ satisfying
\begin{align}
	\bar{B}_j\subseteq \bar{B}_{j+1}, \; j\in \mathbb{N}, \; \lim\limits_{j\to \infty}\bar{B}_j=\mathcal{P}^{n_\zeta}
\end{align}

\begin{algorithm}
	\caption{The Improved Output-based VI Algorithm \cite{Lin2024}} \label{alg2}
	\begin{algorithmic}[1]
		\State Choose $\bar{P}_0=(\bar{P}_0)^T\geq0$. $k,j \gets 0$.
		\Loop
		\State Solve $\bar{H}_k$ from \eqref{vilinearznew}.
		\State
		$\tilde{\bar{P}}_{k+1}\leftarrow \bar{P}_k+\epsilon_k (\bar{H}_k-\bar{P}_kB_\zeta{R}^{-1}B_\zeta^T\bar{P}_k)$
		\If{$\tilde{\bar{P}}_{k+1}\notin\bar{B}_j$}
		\State $\bar{P}_{k+1}\leftarrow \bar{P}_0$. $j\leftarrow j+1$.
		\ElsIf{$||\tilde{\bar{P}}_{k+1}-\bar{P}_k||/\epsilon_k<\varepsilon$}
		\State \Return $\bar{P}_k$ and $-{R}^{-1}B_\zeta^T\bar{P}_k$ as  approximations to ${P}_\zeta^*$  and $K_\zeta^*$, respectively.
		\Else
		\State $\bar{P}_{k+1}\leftarrow \tilde{\bar{P}}_{k+1}$
		\EndIf
		\State $k \leftarrow k+1$
		\EndLoop
	\end{algorithmic}
\end{algorithm}

\begin{rem} A few advantages of Algorithm \ref{alg2} over the one in \cite{Rizvi2019, Rizvi2023} are as follows: First, the one in \cite{Rizvi2019, Rizvi2023} requires $M$ to {\color{black}be of full row rank} while  Algorithm \ref{alg2} does not rely on this assumption. Second, the one in \cite{Rizvi2019, Rizvi2023} needs to solve a sequence of equations with
$\frac{n_\zeta(n_\zeta+1)}{2} + m n_\zeta$ unknown variables. In contrast, by taking advantage of the fact that the matrix $B_{\zeta}$ is known, \cite{Lin2024} reduces the number of unknown variables in (\ref{vilinearznew}) to $\frac{n_\zeta(n_\zeta+1)}{2}$, which is smaller than the number of unknown variables in \cite{Rizvi2019, Rizvi2023} by $m n_\zeta$.
Third, the one in \cite{Rizvi2019, Rizvi2023} requires the initial matrix $\bar{P}_0$ to be positive definite while \cite{Lin2024} allows $\bar{P}_0$ to be positive semi-definite.
	The third advantage is particularly interesting because, in many cases, a Riccati equation only guarantees a positive semi-definite solution.
	Moreover, since the existing state-based PI and VI methods only apply to the case where  $(A,B)$ is stabilizable and $(A,\sqrt{{Q}})$ is observable, they do not directly apply to \eqref{areliz} where the pair $(A_\zeta, B_\zeta)$ is stabilizable but the pair $(A_\zeta, \sqrt{{Q_y}}CM)$ is only detectable. For this reason, they further extended state-based PI and VI methods to the case where  $(A,B)$ is stabilizable and $(A,\sqrt{{Q}})$ is detectable.
\end{rem}

\section{Output Regulation Problem by Parameterized Output-Feedback Control} \label{regulation}

In this section, we first summarize the result of the output regulation problem by the internal model approach based on \cite{Huang2004}. Then we further present a new approach to designing the output-feedback control law to deal with the output regulation problem. This new approach lays the foundation for solving the model-free output regulation problem by output-feedback control law.

\subsection{Existing Results}
Consider the continuous-time linear systems in the following form:

\begin{equation}\label{lisys2}
	\begin{aligned}
		\dot{x}=&Ax+Bu+Ev\\
		y=&Cx\\
		e=&Cx+Fv
	\end{aligned}
\end{equation}
where $x\in \mathbb{R}^n$ is the system state, $u\in\mathbb{R}^m$ is the input, $y\in \mathbb{R}^p$ and $e\in \mathbb{R}^p$ are the measurement output and regulated output, respectively. $v \in \mathbb{R}^{q}$ is generated by the following autonomous system:
\begin{equation} \label{exosys}
	\begin{aligned}
		\dot{v}(t)=Sv(t)
	\end{aligned}
\end{equation}

The output regulation problem is defined as follows:
\begin{pro}\label{StateLORP}
	Given the plant \eqref{lisys2} and the exosystem \eqref{exosys}, design a state-feedback  or output-feedback control law such that the {\color{black}closed-loop system} is exponentially stable with $v$ set to zero and $\lim_{t \to \infty} e(t)=0$.
\end{pro}

The following assumptions are quite standard in solving the linear output regulation problem.
\begin{assmp} \label{ass3}
	$S$ has no eigenvalues with negative real parts.
\end{assmp}

\begin{assmp} \label{ass4}
	For all {\color{black}$\lambda \in \lambda(S)$,}
	\begin{equation}
		\textup{rank}\begin{bmatrix}
			A-\lambda I_n & B\\
			C & \bf{0}
		\end{bmatrix}=n+p.
	\end{equation}
\end{assmp}

The pair $(G_1,G_2)$ is said to be the minimum \textit{p}-copy internal model of $S$ if
\begin{subequations} \nonumber
	\begin{align}
		G_1=\mbox{blockdiag}\underbrace{(\beta,\cdots,\beta)}_{p-tuple},
		G_2=\mbox{blockdiag}\underbrace{(\sigma,\cdots,\sigma)}_{p-tuple},
	\end{align}
\end{subequations}
where $\beta$ is a constant square matrix whose characteristic polynomial equals the minimal polynomial of $S$, and $\sigma$ is a constant column vector such that $(\beta,\sigma)$ is controllable.
A  dynamic compensator of the following form
\begin{align}
		\dot{z}=&G_1z+G_2e \label{im}
	\end{align}
with $z \in \mathbb{R}^{n_z}$
 is called an internal model of {\color{black}the exosystem.} An important implication of the internal model is given by Lemmas 1.26 and  1.27 of \cite{Huang2004},  which is rephrased as follows:

\begin{lem} \label{lem1x}
	Under Assumptions \ref{ass1}, \ref{ass3}, and \ref{ass4}, let $(G_1,G_2)$ be the minimal \textit{p}-copy internal model of $S$. Then, \\
	(i) The pair  $\begin{bmatrix}
		A & \bf{0}\\
		G_2C & G_1
	\end{bmatrix}$ and $\begin{bmatrix}
		B\\
		\bf{0}
	\end{bmatrix}$ is stabilizable. \\
	(ii) For any $V \in \mathbb{R}^{p \times q}$, if the following matrix equation  \\
	\begin{equation}\label{sylv2}
		Z S = G_1 Z + G_2 V
	\end{equation}
	has a solution  $Z \in \mathbb{R}^{n_z \times q}$, then {\color{black}$V = \mathbf{0}$.}
	
\end{lem}

Let $\xi=\col (x,z)$.  Combining \eqref{lisys2} with \eqref{im} gives the following so-called augmented system:
\begin{subequations}\label{aug}
	\begin{align}
		\dot{\xi} &=\begin{bmatrix}
			A & \bf{0}\\
			G_2C & G_1
		\end{bmatrix}\xi+
		\begin{bmatrix}
			B\\
			\bf{0}
		\end{bmatrix}u+
		\begin{bmatrix}
			E\\
			G_2F
		\end{bmatrix}v\\
		e&=\begin{bmatrix}
			C & \bf{0}
		\end{bmatrix}\xi+Fv
	\end{align}
\end{subequations}
Let $Y=\begin{bmatrix}
   	A& \bf{0}\\ {G}_2C&{G}_1
   \end{bmatrix}$,  {\color{black}$J=\begin{bmatrix}
   	B\\ \bf{0}
   \end{bmatrix}$}, and $\bar{E} = \begin{bmatrix}
   		E\\G_2F
   	\end{bmatrix}$. Then \eqref{aug} can be put into the following compact form:
   \begin{align}\label{augcom}
   	\dot{\xi}=&Y\xi+Ju+ \bar{E}v
   \end{align}

Then, based on Lemma 1.26 and Remark 1.28 of \cite{Huang2004}, the solvability of Problem \ref{StateLORP} is summarized by the following theorem.

\begin{thm} \label{thm1}
	Under Assumptions \ref{ass1}, \ref{ass3}, and \ref{ass4}, let $(G_1,G_2)$ be the minimal \textit{p}-copy internal model of $S$.  Let  $K=\begin{bmatrix}
		K_x & K_z
	\end{bmatrix}$ be such that $Y+ JK$ is Hurwitz. Then Problem \ref{StateLORP} is solved by the following dynamic state-feedback control law
	\begin{subequations} \label{uz}
	\begin{align}
		u=&K_xx+K_zz \label{u}\\
		\dot{z}=&G_1z+G_2e \label{z}
	\end{align}
\end{subequations}
Moreover, let $L$ be such that $A- LC$ is Hurwitz.  Then Problem \ref{StateLORP} is solved by the following dynamic output-feedback control law
		\begin{subequations} \label{uzo}
		\begin{align}
			u &=K_x \hat{x}+K_zz \\
			\dot{\hat{x}} & = A \hat{x} + B u + L (y - C \hat{x}) \label{uzob} \\
			\dot{z}& =G_1z+G_2e
		\end{align}
	\end{subequations}
\end{thm}
\begin{rem}
  The state-feedback control law \eqref{uz} only involves one design parameter $K$.  For any $Q_{\xi}=Q_{\xi}^T\geq 0, R=R^T>0$ with $(Y,\sqrt{Q}_{\xi})$ observable, one
  can define a Riccati equation of the following form
 \begin{align}\label{riccati}
	Y^TP_{\xi}^*+P_{\xi}^*Y+Q_{\xi}-P_{\xi}^*JR^{-1}J^TP_{\xi}^*=0
\end{align}
Since $(Y, J)$ is stabilizable,  \eqref{riccati} has a unique positive definite solution $P_{\xi}^*$.
Then a stabilizing feedback gain is given by $K_{\xi}^*=-R^{-1}J^TP_{\xi}^*$. Based on \eqref{riccati} and \eqref{augcom}, one can further develop a data-driven approach
to obtain $K$ without knowing the parameters of the plant \cite{LH2024}.

\end{rem}	

\begin{rem}\label{remop2}
Since,  under Assumption \ref{ass3}, the spectra of $S$ and $Y+JK $ are disjoint,
	by Proposition A.2 in \cite{Huang2004},   the following Sylvester equation admits a unique solution $X_c \in  \mathbb{R}^{(n +n_z) \times q}$:
\begin{align}\label{Xc}
		X_c S=&(Y+J K)X_c+\bar{E}
		\end{align}
Partition $X_c$ as
	  $X_c   = \left[
	\begin{array}{c}
		X   \\
		{Z}
	\end{array}
	\right]$  with $X   \in \mathbb{R}^{n \times q} $, and ${Z}   \in \mathbb{R}^{n_z \times q}$. Then \eqref{Xc} implies

		\begin{align}
			{Z}    S=& {G}_1 {Z}     + {G}_2 V  \label{ire1oc0}
		\end{align}
	where
	\begin{equation*} %\label{ire2o}
		V   = {C} X     + {F}
	\end{equation*}
	Since  equation \eqref{ire1oc0}
	is in the form (\ref{sylv2}), by
	Lemma \ref{lem1x}, {\color{black}$V=\mathbf{0}$, that is,
\begin{align}
	\mathbf{0}=&\begin{bmatrix} C & \mathbf{0} \end{bmatrix} X_c+F \label{xz}
\end{align}}

Define $\bar{\xi}=\xi-X_cv$, $\bar{u}=u-KX_cv$.  By \eqref{xz} and \eqref{Xc}, we have
\begin{align}\label{barxic}
	\dot{\bar{\xi}}= Y \bar{\xi} +J\bar{u}
\end{align}
Thus, the control law as given by $\bar{u}^*=K_{\xi}^* \bar{\xi}$ can be interpreted as the  solution to the following LQR problem:
\begin{pro}\label{StateLQR}
\begin{equation}\label{op}
		\begin{aligned}
			& \min_{\bar{u}}\int_0^\infty(|\bar{\xi}|_{{Q}_{\xi}}+|\bar{u}|_R)dt, {Q}_{\xi}\geq0, R>0\\
			& \text{subject to}\quad \eqref{barxic}
		\end{aligned}
	\end{equation}
\end{pro}
\end{rem}

\subsection{Output Regulation Problem by Parameterized Output-Feedback Control}

{\color{black} Since the observer \eqref{uzob} in the output-feedback control law \eqref{uzo} is defined by unknown matrices $A$, $B$, $C$, and $L$, the control law \eqref{uzo} cannot be implemented. In this subsection, we will replace  \eqref{uzob}  by a known  observer given in \eqref{dyz}.  In conjunction with the VI based learning approach, this control law will lead to a model-free approach to solving the output regulation problem by an output-feedback control law.}

Let  $K=\begin{bmatrix}
	K_x & K_z
\end{bmatrix}$ be such that $Y+JK$ is Hurwitz and $L$ be such that $A-LC$ is Hurwitz.  We consider the following dynamic output-feedback control law
\begin{subequations} \label{uz2x}
\begin{align}
		u=&K_x M \zeta+K_z z \\
		\dot{\zeta}=&(I_{m+p}\otimes \mathcal{A}) \zeta+ B_\zeta u+E_\zeta y  \label{uz2xb}\\
		\dot{z}=&G_1z+G_2e
	\end{align}
\end{subequations}

{\color{black}Let  $x_c=\mbox{col} (x, \zeta, z)$. Combining \eqref{lisys2} with \eqref{uz2x} gives the following closed-loop system:
\begin{subequations}\label{xc}
	\begin{align}
		\dot{x}_c&= A_c x_c + B_c v \\
		e&= C_c x_c+ F v
	\end{align}
\end{subequations}
where
$A_c =\begin{bmatrix}
	A  & B K_x M & B K_z \\
	E_\zeta C  &	I_{m+p}\otimes \mathcal{A} +B_\zeta K_x M  & B_\zeta K_z\\
	G_2C & \bf{0}  & G_1
\end{bmatrix}$, $W_c=\mbox{blockdiag}(I_n, M, I_{n_z})$, $B_c =\begin{bmatrix}
	E\\
	\bf{0} \\
	G_2 F
\end{bmatrix}$, and $C_c=\begin{bmatrix}
	C & \bf{0} & \bf{0}
\end{bmatrix}$.}

\begin{lem} \label{lem3}
	Under Assumptions \ref{ass1}-\ref{ass4}, there exist $K_x$, $K_z$ and $L$ such that $A_c$ is Hurwitz.
\end{lem}

\begin{Prf}
	By Lemma \ref{lem1x}, under Assumptions \ref{ass1}, \ref{ass3}, and \ref{ass4}, the pair $(Y, J)$ is stabilizable. Thus, there exists a $K = \begin{bmatrix}K_x & K_z \end{bmatrix}$ such that $Y + J K = \begin{bmatrix}
		A  + B K_x & B K_z  \\
		G_2C & G_1
	\end{bmatrix} $ is Hurwitz, and under Assumptions \ref{ass2}, there exists $L$  such that $A-LC$ is Hurwitz.
	Let ${W}_c = \mbox{blockdiag}(I_{n}, M, I_{n_z})$  and
		\begin{equation}\label{barac}
		 {A}_\eta = \begin{bmatrix}
			A  & B K_x  & B K_z \\
			L C  &	A - L C +B K_x  & B K_z\\
			G_2C & \mathbf{0}  & G_1
		\end{bmatrix}
		\end{equation}
Using the following identities from  Lemma 5 of \cite{Lin2024}
\begin{align}\label{rela1}
	M (I_{m+p}\otimes \mathcal{A}) =&(A-LC)M\\ \label{rela2}
	M\begin{bmatrix}
		I_m\otimes b\\ {\bf0}
	\end{bmatrix}=&B\\ \label{rela3}
	M\begin{bmatrix}
		{\bf0}\\I_p\otimes b
	\end{bmatrix}=&L
\end{align}
gives
	\begin{equation}
		W_c {A}_c = \begin{bmatrix}
			A  & B K_x  & B K_z \\
			L C  &	A - L C +B K_x  & B K_z\\
			G_2C & \mathbf{0}  & G_1
		\end{bmatrix}W_c = A_\eta W_c
\end{equation}
Thus, letting $\eta= {W}_c x_c = \mbox{col} (x, M\zeta, z)$	gives	
\begin{equation}\label{barxc}
		\begin{aligned}
			\dot{\eta}&= W_c {A}_c x_c + W_c B_c v \\
			&= A_\eta 	\eta+ W_c B_c v
		\end{aligned}
\end{equation}

Let
	\begin{equation}
		T = \begin{bmatrix}
			I_n   &  \mathbf{0}  &  \mathbf{0} \\
		- I_n &	I_n  & \mathbf{0}\\
		\mathbf{0} & \mathbf{0}  & I_{n_z}
		\end{bmatrix}
	\end{equation}
	Then
	\begin{equation}
		T {A}_\eta T^{-1}  = \begin{bmatrix}
			A  + B K_x & B K_x  & B K_z \\
		\mathbf{0}   &	A - L C   &  \mathbf{0} \\
			G_2C & \mathbf{0}  & G_1
		\end{bmatrix}
	\end{equation}
	
	Since  $\begin{bmatrix}
		A  + B K_x & B K_z  \\
		G_2C & G_1
	\end{bmatrix}$ and $A-LC$ are Hurwitz, ${A}_\eta$ is Hurwitz.
	Thus, $\eta = \mbox{col} (x, M\zeta, z)$ converges to the origin exponentially when $v$ set to zero.

	By \eqref{uz2x}, with $v$ set to zero,  $\zeta$ satisfies
	\begin{equation}
		\begin{aligned}
			\dot{\zeta} =&(I_{m+p}\otimes \mathcal{A})\zeta+ B_\zeta u  + E_\zeta y \\
			=&(I_{m+p}\otimes \mathcal{A})\zeta+ B_\zeta K_x M\zeta+B_\zeta K_z z+E_\zeta C x
		\end{aligned}
	\end{equation}
	Since $(I_{m+p}\otimes \mathcal{A})$ is Hurwitz, and $\eta$  decays to zero exponentially when $v$ set to zero, $\zeta$ decays to zero exponentially.  Thus, $x_c$
	decays to zero exponentially when $v$ set to zero,  which implies  $A_c$ is Hurwitz.
\end{Prf}

Now we are ready to establish the following result.
\begin{thm}\label{thm2}
Under Assumptions \ref{ass1}-\ref{ass4}, let  $K = [K_x, K_z]$ and $L$ be such that $Y+JK$ and  $A- LC$  Hurwitz. Then  the  dynamic output-feedback control law \eqref{uz2x} {\color{black}
solves Problem \ref{StateLORP}.}
\end{thm}

\begin{Prf}
	By Lemma \ref{lem3}, there exists $K_x$, $K_z$ and $L$  such that $A_c $ is Hurwitz.
	Under Assumption \ref{ass3}, the spectra of $S$ and $A_{c} $ are disjoint. Thus,
	by Proposition A.2 in \cite{Huang2004},   the following Sylvester equation has a unique solution $\tilde{X}_c$:
	\begin{align}
		\tilde{X}_c   S &= {A}_c  \tilde{X}_c   + B_c   \label{ire1cc1}
	\end{align}
Partition $\tilde{X}_c$ as  $\left[
	\begin{array}{c}
		\tilde{X}   \\
		\hat{Z}   \\
		\tilde{Z}
	\end{array}
	\right]$  with $\tilde{X}   \in \mathbb{R}^{n \times q} $, $\hat{Z}  \in \mathbb{R}^{n (p+m) \times q}$, and $\tilde{Z}   \in \mathbb{R}^{n_z \times q}$.
	Then the Sylvester equation  \eqref{ire1cc1} can be expanded as follows:
	\begin{subequations}\label{ire1o}
		\begin{align}
			\tilde{X}   S=& {A}   \tilde{X}   + B(  K_x M \hat{Z}   + K_z \tilde{Z}   )   + {E} \label{ire1oa}\\
			\hat{Z}    S=&  	E_\zeta C  \tilde{X} +	(I_{m+p}\otimes \mathcal{A} +B_\zeta K_x M ) \hat{Z}    +  B_\zeta K_z \tilde{Z}   \label{ire1ob}\\
			\tilde{Z}    S=& {G}_1 \tilde{Z}     + {G}_2 \tilde{V}   \label{ire1oc}
		\end{align}
	\end{subequations}
	where
	\begin{equation*} %\label{ire2o}
		\tilde{V}   = {C} \tilde{X}     + {F}
	\end{equation*}
	Since  equation \eqref{ire1oc} is in the form (\ref{sylv2}), by
	Lemma \ref{lem1x}, $\tilde{V}  = \mathbf{0}$, which implies
	\begin{align}
	\mathbf{0}&= C_c \tilde{X}_c   + F. \label{ire1cc2}
	\end{align}
	Let $\bar{x}_c = x_c - \tilde{X}_c v$. Then using  \eqref{ire1cc1} and \eqref{ire1cc2} gives the following:
	\begin{subequations}\label{ire1ox}
		\begin{align}
			\dot{\bar{x}}_c  =& A_c \bar{x}_c  \\
			e =&  C_c 	\bar{x}_c
		\end{align}
	\end{subequations}
Since 	$A_c $ is Hurwitz, $\bar{x}_c$ and  hence  $e$ tend to zero exponentially. The proof is complete.
\end{Prf}

\section{Main Results}\label{main results}

In this section, we will further consider designing a model-free output-based control law for solving the output regulation problem as formulated below:

\begin{pro}\label{LORP}
Given the plant \eqref{lisys2} and the exosystem \eqref{exosys}, design a dynamic output-feedback control law such that the  {\color{black} closed-loop system} is exponentially stable with $v$ set to zero and $\lim_{t \to \infty} e(t)=0$ without knowing the plant parameters $A$, $B$, $C$, $E$, and $F$.
\end{pro}

{\color{black}Since we assume the  internal model \eqref{im} is known,  we need one more assumption as follows:

\begin{assmp} \label{ass5} The minimal polynomial of $S$ is known.
\end{assmp}
}

\subsection{Output-based Output Regulation Problem by Modified Internal Model Approach}\label{sec3-1}

{\color{black}Since we deal with the asymptotic tracking and disturbance rejection problems by a control law independent of $v$,  the methods proposed in \cite{Chen2022} and \cite{Xie2023} do not apply to our problem.  We need to further empower two new features in our approach. First, note that,  due to the presence of the persistent exogenous signal $v$, the state $\hat{x}$ of the Luenberger observer \eqref{observer} does not converge to the state $x$ of the plant. In fact, the following relation holds:
\begin{align}
\dot{x}-\dot{\hat{x}}=(A-LC)(x-\hat{x})+Ev
\end{align}
Thus, $M\zeta$ does not converge to the state $x$ of the plant, either. Therefore, we need to find out the relation between $x$ and $\zeta$.
The second one is more subtle. For the output regulation problem, we need to design an output-based control law to stabilize the augmented system \eqref{aug} instead of the original plant. Thus, we need to derive a so-called augmented auxiliary  system for the augmented system \eqref{aug}, which plays the same role as the auxiliary  system \eqref{dyz2} does to the LQR problem for \eqref{lisys}.

Let us first find out the relation between $x$ and $\zeta$.}
Under Assumption \ref{ass2}, there exists a $L$ such that $A-LC$ is Hurwitz. Then,  by Proposition A.2 in \cite{Huang2004},  the following Sylvester equation
\begin{align}\label{SylX}
X'S=(A-LC)X'+E
\end{align}
admits a unique solution $X' \in \mathbb{R}^{n \times q}$ since the spectra of $S$ and $A-LC$ are disjoint.
Let $e_x(t)=M\zeta(t)+X'v(t)-x(t)$. {\color{black}Then by \eqref{rela1}-\eqref{rela3} and \eqref{SylX}, we have
\begin{align} \label{delx}
\dot{e}_x=&M\dot{\zeta}+X'\dot{v}-\dot{x}\notag \\
=&M((I_{m+p}\otimes \mathcal{A})
\zeta+\begin{bmatrix}
I_m\otimes b\\ {\bf0}
\end{bmatrix}u+\begin{bmatrix}
{\bf0}\\I_p\otimes b
\end{bmatrix}y)\notag \\
&+X'Sv-Ax-Bu-Ev \notag \\
=&(A-LC)M\zeta+Bu+Ly+X'Sv-Ax-Bu-Ev \notag \\
=&(A-LC)(M\zeta-x)+(X'S-E)v \notag \\
=&(A-LC)(M\zeta-x)+(A-LC)X'v \notag \\
=&(A-LC)e_x
\end{align}}
Thus,  \eqref{delx} is exponentially stable.

Thus, we have the following relationship:
\begin{align} \label{x4}
x(t)=M\zeta(t)+X'v(t)-e_x(t)
\end{align}
where $e_x(t)$ converges to zero exponentially.

\begin{rem}
Reference  \cite{Chen2022} only considered the tracking problem where $E=0$. In this case, $\hat{x}$ will converge to x.
{\color{black}On the other hand,
reference \cite{Xie2023} constructed a known observer based on the following Luenberger observer:
\begin{align}\label{observer3}
\dot{\hat{x}}= &A\hat{x}+Bu+Ev+L(y-C\hat{x})
\end{align}
which depends on $v$. As a result, the control law in \cite{Xie2023} also depends on $v$ and has high dimension. More detailed comparisons will be given in Remarks \ref{rem6} and \ref{rem8}.}
\end{rem}

Then we need to define an augmented auxiliary  system so that the desired control gain can be obtained by solving a Riccati equation associated with this augmented auxiliary  system. For this purpose,
substituting  $y=Cx=CM\zeta+CX'v-Ce_x$ into \eqref{dyz} yields the following dynamics:
\begin{align}\label{newzeta}
\dot{\zeta}=A_\zeta\zeta+B_\zeta u+E_\zeta CX' v-E_\zeta Ce_x
\end{align}
Also, substituting \eqref{x4} into \eqref{z} yields the following:
\begin{equation}\label{z3}
\begin{aligned}
\dot{z}=G_2CM\zeta+G_1z+G_2(CX'+F)v-G_2Ce_x
\end{aligned}
\end{equation}

Let $\rho=\mbox{col} (\zeta, z)$. Then, the compact form of   \eqref{newzeta} and \eqref{z3} is as follows:
\begin{equation}\label{rhoz2}
\begin{aligned}
\dot{\rho}=A_\rho \rho+B_\rho u+E_\rho v+D_{\rho}e_x
\end{aligned}
\end{equation}
where $\rho \in \mathbb{R}^{n_\rho}$ with $n_{\rho}=n(m+p)+n_z$,
\begin{align}
A_\rho=&\begin{bmatrix}
A_\zeta & \bf{0}\\
G_2CM & G_1
\end{bmatrix},
B_\rho=\begin{bmatrix}
B_\zeta \\
\bf{0}
\end{bmatrix} \notag\\
E_\rho=&\begin{bmatrix}
E_\zeta CX' \\
G_2(CX'+F)
\end{bmatrix},
D_{\rho}=\begin{bmatrix}
-E_\zeta C\\
-G_2C
\end{bmatrix} \notag
\end{align}

It is interesting to note that while $A_\rho$ is unknown, $B_\rho$ is known. In what follows, system \eqref{rhoz2} will play the same role in solving the output regulation problem for \eqref{lisys2} as what \eqref{dyz2bf} has
played in solving the LQR problem of \eqref{lisys}. For this purpose, we first establish the following result.
\begin{thm}\label{thm3}
Under Assumptions \ref{ass1}-\ref{ass5}, \\
(i) The pair $(A_\rho, B_\rho)$ is stabilizable.\\
(ii) Let $K_\rho$ be any matrix such that $A_\rho + B_\rho K_\rho$ is Hurwitz and $L$ be such that $A- LC$ is Hurwitz.
Then,  Problem \ref{LORP} is solved by the following dynamic output-feedback control law:
\begin{equation}\label{uz2}
\begin{aligned}
u=&K_\rho \rho \\
\dot{\zeta}=&(I_{m+p}\otimes \mathcal{A}) \zeta+ B_\zeta u+E_\zeta y \\
\dot{z}=&G_1z+G_2e
\end{aligned}
\end{equation}
\end{thm}

\begin{Prf}
Part (i). Under Assumptions \ref{ass1} to \ref{ass5}, let  $\begin{bmatrix}K_x & K_z \end{bmatrix}$ be such that $\begin{bmatrix}
A & \bf{0}\\
G_2C & G_1
\end{bmatrix}+\begin{bmatrix}
B\\
\bf{0}
\end{bmatrix}\begin{bmatrix}K_x & K_z \end{bmatrix}$ is Hurwitz and $L$ be such that $A- LC$ is Hurwitz.

 Let $K_\rho = \begin{bmatrix}K_x & K_z \end{bmatrix}W$, where $W=\mbox{blockdiag}(M, I_{n_z})$. With $K_\rho$ designed this way, control law  \eqref{uz2}  is the same as control law \eqref{uz2x}. Thus,  by Lemma \ref{lem3},  $x$, $\zeta$, $z$  all converge to the origin exponentially when $v$ is set to zero.

Now substituting $u=K_{\rho} {\rho}$  into \eqref{rhoz2} and setting $v=0$ gives
\begin{equation}\label{rhou}
\begin{aligned}
\begin{bmatrix}
\dot{\rho}\\
\dot{e}_x
\end{bmatrix}=&\begin{bmatrix}
A_\rho+B_\rho K_\rho & D_\rho\\
\bf{0} & A- LC
\end{bmatrix}\begin{bmatrix}
\rho\\
e_x
\end{bmatrix}
\end{aligned}
\end{equation}
Since we have shown that  $\zeta$, $z$ and $e_x$ all converge to the origin exponentially when $v$ is set to zero, the matrix $A_\rho+B_\rho K_\rho =
\begin{bmatrix}
A_\zeta+B_\zeta K_xM & B_\zeta K_z\\
G_2CM & G_1
\end{bmatrix}$
must be Hurwitz.   Thus, the pair $(A_\rho, B_\rho)$ must be stabilizable.

Part (ii). Since the pair $(A_\rho, B_\rho)$ is stabilizable, there exists  $K_\rho$  such that $A_\rho + B_\rho K_\rho$ is Hurwitz.
Thus, by setting $v=0$, $\col (\rho, e_x)$ satisfies \eqref{rhou}.
Since, $\zeta$, $z$ and $e_x$ are all converge to the origin exponentially,
by \eqref{x4},  when $v$ is set to $0$, $x$ becomes
\begin{align}
x(t)=M\zeta(t)-e_x(t)
\end{align}
Thus, $x$ also decays to $0$ exponentially.

Let  $K_\rho = [K_{1 \rho}, K_{2 \rho}]$ with  $K_{2 \rho} \in \mathbb{R}^{n_z}$ and $x_c=\mbox{col} (x, \zeta, z)$,   and denote the closed-loop system  as follows:
\begin{subequations}\label{xcxy}
	\begin{align}
		\dot{x}_c&= \bar{A_c} x_c + B_c v \\
		e&= C_c x_c+ F v
	\end{align}
\end{subequations}
where $\bar{A_c}=\begin{bmatrix}
	A  & B K_{1 \rho} & B K_{2 \rho} \\
	E_\zeta C  &	I_{m+p}\otimes \mathcal{A} +B_\zeta K_{1 \rho}  & B_\zeta K_{2 \rho}\\
	G_2C & \bf{0}  & G_1
\end{bmatrix}$,  and $B_c$ and $C_c$ are as defined in  \eqref{xc}.
Since we have shown that $x_c$ tends to zero exponentially when $v$ is set to zero, $\bar{A_c}$ must be Hurwitz.
The rest of the proof is quite similar to the proof of
Theorem \ref{thm2} and will be only sketched below. In fact, since $\bar{A_c}$ is Hurwitz,
the following Sylvester equation has a unique solution $\hat{X}_c$:
	\begin{align}
		\hat{X}_c   S &= \bar{A}_c \hat{X}_c   + B_c   \label{ire1cc12}
	\end{align}
Again, partition $\hat{X}_c   = \left[
	\begin{array}{c}
		\hat{X}   \\
		\hat{Z}'   \\
		\tilde{Z}'
	\end{array}
	\right]$  with $\hat{X}   \in \mathbb{R}^{n \times q} $, $\hat{Z}'  \in \mathbb{R}^{n (p+m) \times q}$, and $\tilde{Z}'   \in \mathbb{R}^{n_z \times q}$.
	 Then the Sylvester equation  \eqref{ire1cc12}  implies the following
\begin{subequations}\label{ire1oxy}
		\begin{align}
			\hat{X} S=& {A}   \hat{X}   + B( K_{1 \rho} \hat{Z}'   + K_{2 \rho}  \tilde{Z}'   )   + {E} \label{ire1oaxy}\\
			\hat{Z}'  S=&  	E_\zeta C  \hat{X} +	(I_{m+p}\otimes \mathcal{A} +B_\zeta  K_{1 \rho}) \hat{Z}'    +  B_\zeta K_{2 \rho} \tilde{Z}'   \label{ire1obxy}\\
			\tilde{Z}' S=& {G}_1 \tilde{Z}'     + {G}_2 V'   \label{ire1ocxy}
		\end{align}
	\end{subequations}
	where
	\begin{equation*} %\label{ire2o}
		V'   = {C}_c  \hat{X}_c     + {F}
	\end{equation*}
	Since  equation \eqref{ire1ocxy}
	is in the form (\ref{sylv2}), by
	Lemma \ref{lem1x}, $V' = \bf{0}$, which implies
	\begin{align}
		\bf{0}&= C_c    \hat{X}_c   + F. \label{ire1cc2xy}
	\end{align}
	Let $\tilde{x}_c = x_c - \hat{X}_c v$. Then using  \eqref{ire1cc12} and \eqref{ire1cc2xy} gives the following:
	\begin{subequations}\label{ire1oxtil}
		\begin{align}
			\dot{\tilde{x}}_c  =& \bar{A}_c \tilde{x}_c  \\
			e =&  C_c 	\tilde{x}_c
		\end{align}
	\end{subequations}
Since 	$\bar{A}_c $ is Hurwitz, $\tilde{x}_c$ and  hence  $e$ tends to zero exponentially. The proof is complete.
\end{Prf}

\begin{rem}

Since $e_x$ decays to 0 exponentially, the term $D_\rho e_x$ in \eqref{rhoz2} can be ignored  after some finite $t_0>0$. Thus, when $t\geq t_0$,
\eqref{rhoz2} is simplified to the following
\begin{equation}\label{rhoz5}
\begin{aligned}
\dot{\rho}=A_\rho \rho+B_\rho u+E_\rho v
\end{aligned}
\end{equation}
We call \eqref{rhoz5} the augmented auxiliary  system.

What makes Theorem \ref{thm3} more interesting than Theorem \ref{thm2} is that it does not require $K_\rho = [K_x M, K_z]$ which is not realizable since $M$ is unknown. Thus it reduces Problem \ref{LORP} to the stabilization of the  augmented auxiliary  system \eqref{rhoz5}
so that
the desired control gain can be obtained by solving a Riccati equation associated with the augmented auxiliary  system.
In fact, let $Q_\rho \geq 0$ be such that $(A_\rho, \sqrt{Q_\rho})$ is  observable and  define  an algebraic Riccati equation as follows:
\begin{equation}\label{AREMrho}
\begin{aligned}
A_\rho^TP_\rho^{*}+P_\rho^{*}A_\rho+Q_\rho-P_\rho^{*}B_\rho R^{-1}B_\rho^TP_\rho^{*}=0.
\end{aligned}
\end{equation}
Since $(A_\rho, B_\rho)$ is stabilizable and $(A_\rho, \sqrt{Q_\rho})$ is observable, \eqref{AREMrho} admits a unique positive definite solution $P_\rho^{*} $ \cite{kucera1972, Wonham}. Based on the solution $P_\rho^{*} $ to \eqref{AREMrho}, it is possible to obtain a stabilizing gain $K_\rho^{*} $ without knowing the parameters of the plant as will be detailed in the next subsection.

\end{rem}

\begin{rem}\label{remop3x}
Similar to Remark \ref{remop2}, we can also relate the Ricatti equation \eqref{AREMrho} to some LQR problem defined as follows. In fact, since,  under Assumption \ref{ass3}, the spectra of $S$ and $A_\rho + B_\rho  K_\rho $ are disjoint,
	by Proposition A.2 in \cite{Huang2004},
	 the following Sylvester   equation admits a unique solution $X_\rho \in  \mathbb{R}^{n_\rho \times q}$:
\begin{align}\label{Xrho}
		X_\rho S=&(A_\rho + B_\rho  K_\rho)X_\rho+ {E}_\rho
		\end{align}
Partition $X_\rho$ as
	  $X_\rho   = \left[
	\begin{array}{c}
		X_\zeta   \\
		X_z
	\end{array}
	\right]$  with $X_\zeta   \in \mathbb{R}^{n_\zeta \times q} $, and $X_z   \in \mathbb{R}^{n_z \times q}$. Then \eqref{Xrho} implies
\begin{equation}
\begin{aligned}
X_z S=& G_1 X_z + G_2 (\begin{bmatrix} CM & {\color{black}\bf{0}} \end{bmatrix} X_\rho + C X' +{F} ) \label{ire1oc0x}
\end{aligned}
\end{equation}
Since  equation \eqref{ire1oc0x} is in the form (\ref{sylv2}), by Lemma \ref{lem1x},
\begin{equation}
\begin{aligned}\label{xz3x}
\bf{0}=&\begin{bmatrix} CM & \bf{0} \end{bmatrix}X_\rho + C X' +{F}
\end{aligned}
\end{equation}
Let $\bar{\rho} =\rho-X_\rho v$, $\bar{u}=u-K_\rho X_\rho v$.  By \eqref{xz3x} and \eqref{Xrho}, we have
\begin{align}\label{barrho}
	\dot{\bar{\rho}}= A_\rho \bar{\rho} +B_\rho  \bar{u}
\end{align}
Thus, the control law as given by $\bar{u}^*=K^*_\rho  \bar{\rho}$ can be interpreted as the  solution to the following LQR problem:

\begin{pro}\label{OPFBLQR}
\begin{equation}\label{op3}
		\begin{aligned}
			& \min_{\bar{u}}\int_0^\infty(|\bar{\rho}|_{{Q}_\rho}+|\bar{u}|_R)dt, {Q}_\rho \geq0, R>0\\
			& \text{subject to}\quad \eqref{barrho}
		\end{aligned}
	\end{equation}

\end{pro}
\end{rem}

Although the dynamics of $\rho$ is unknown, $\rho$ is available. Thus, it is possible  to obtain a stabilizing feedback gain $K^*_\rho$ by solving the algebraic Riccati equation \eqref{AREMrho} with the VI based method.

\subsection{VI Method for Solving \eqref{AREMrho}}\label{ADP}

To obtain the VI based algorithm for solving \eqref{AREMrho}, let $H_k^{\rho}=A_\rho^TP_{k}^{\rho}+P_{k}^{\rho} A_\rho$.
Integrating $\frac{ d (\rho^T (t) {{P}^\rho_k} \rho (t))}{dt}$ along the solution of  \eqref{rhoz5}
 gives
\begin{equation}\label{RZIRL}
\begin{aligned}
&|\rho(t+\delta t)|_{P_{k}^{\rho}}-|\rho(t)|_{P_{k}^{\rho}}\\
=&\int_{t}^{t+\delta t} \left( |\rho|_{H_k^{\rho}}  + 2u^TB_\rho^TP_{k}^{\rho}\rho +2v^TE_\rho^TP_{k}^{\rho} \rho \right) d\tau
\end{aligned}
\end{equation}

For any vectors $a$, $b$ and any integer $s\in \mathbb{N}_+$, let
\begin{align}
\Gamma_{ab}=&[\int_{t_0}^{t_1}a\otimes b d\tau , \int_{t_1}^{t_2}a\otimes b d\tau, \cdots , \int_{t_{s-1}}^{t_s}a\otimes b d\tau]^T
\end{align}

Then \eqref{RZIRL} implies
\begin{equation}\label{linear equationR}
\begin{aligned}
\Psi^{\rho}\begin{bmatrix}
\text{vecs}(H_k^{\rho})\\
\text{vec}(E_\rho^TP_{k}^{\rho})
\end{bmatrix}=\Phi_{k}^{\rho}
\end{aligned}
\end{equation}
where
\begin{equation}
\begin{aligned}
\Psi^{\rho}&=[I_{\rho \rho}, 2\Gamma_{\rho v}]\\
\Phi_{k}^{\rho}&=\delta_{\rho}\text{vecs}(P_{k}^{\rho})-2\Gamma_{\rho (B_\rho u)}\text{vec}(P_{k}^{\rho})
\end{aligned}
\end{equation}

The solvability of \eqref{linear equationR} is guaranteed by the following result.
\begin{lem} \label{rankVI3}
The matrix $\Psi^{\rho}$ has full column rank if
\begin{align}\label{rank3}
\textup{rank}([I_{\rho \rho}, \Gamma_{\rho v}])=& \frac{n_{\rho}(n_{\rho}+1)}{2}+qn_{\rho}
\end{align}
\end{lem}

The VI algorithm for  approximating $-R^{-1}B_\rho^T P_\rho^{*}$ is summarized in Algorithm \ref{vialg3}, where $\epsilon_k$ and $\varepsilon$ are defined in the same way as those in Section \ref{stateLQR},  $\{\tilde{B}_j\}_{q=0}^{\infty}$ is a collection of bounded subsets in $\mathcal{P}^{n_\rho}$ satisfying
\begin{align}
\tilde{B}_j\subseteq \tilde{B}_{j+1}, \; j\in \mathbb{N}, \; \lim\limits_{j\to \infty}\tilde{B}_j=\mathcal{P}^{n_\rho}
\end{align}

\begin{algorithm}
\caption{Model-free VI Algorithm for Solving \eqref{AREMrho}} \label{vialg3}
\begin{algorithmic}[1]
\State Apply any locally essentially bounded initial input $u^0$.
Choose a proper $t_0 > 0$ such that $\|e_x\|$ is small enough.
Collect data from $t_0$ until the rank condition \eqref{rank3} is satisfied.
\State Choose $P_{0}^{\rho}=(P_{0}^{\rho})^T>0$. $k,j \gets 0$
\Loop
\State Substitute $P_{k}^{\rho}$ into \eqref{linear equationR} to solve $H_k^{\rho}$ and $K_{k}^{\rho}=-R^{-1}B_\rho^TP_{k}^{\rho}$.
\State \begin{equation}\label{A3}
\begin{aligned}
\tilde{P}_{k+1}^{\rho} \gets P_{k}^{\rho}+\epsilon_k(H_k^{\rho}+Q_\rho-(K_{k}^{\rho})^TRK_{k}^{\rho})\\
\end{aligned}
\end{equation}
\If{$\tilde{P}_{k+1}^{\rho} \notin \tilde{B}_j $}
    \State $P_{k+1}^{\rho} \gets P_{0}^{\rho}$. $j \gets j+1$.
\ElsIf{$||\tilde{P}_{k+1}^{\rho}-P_{k}^{\rho}||/\epsilon_k<\varepsilon$}
    \State \textbf{return} $K_{k}^{\rho}$ as an approximation to $-R^{-1}B_\rho^TP_\rho^{*}$
\Else
	\State $P_{k+1}^{\rho} \gets \tilde{P}_{k+1}^{\rho}$
\EndIf
\State $k \gets k+1$
\EndLoop
\end{algorithmic}
\end{algorithm}

\begin{rem}
In deriving \eqref{linear equationR}, we have taken advantage
 of the fact that  $B_\rho$ is known. Thus,  $K_k^{\rho} =  -R^{-1} B_\rho^TP_{k}^{\rho} $ is also known. As a result, we don't need to solve $K_k^{\rho}$ at each iteration as opposed to other approaches such as
that of \cite{Chen2022}.
\end{rem}

{\color{black}\begin{rem} \label{remgao}
It is noted that the quantity $\Gamma_{\rho v}$ depends on the exogenous signal $v$. Nevertheless, under Assumption \ref{ass5},  we can always assume $v$ is known.
In fact, as pointed out in \cite{Gao2016}, under Assumption \ref{ass5}, there exist a known vector $\hat{v} (t)$, a known matrix $\hat{S}$ whose characteristic polynomial is equal to the minimal polynomial of $S$, and an unknown vector $\hat{C}$ depending on $v (0)$ such that,
\begin{align}
\dot{\hat{v}} = \hat{S} \hat{v},~~ v = \hat{C} \hat{v}
\end{align}
As a result, equation \eqref{lisys2} can be recast as follows:
\begin{equation}\label{lisys2re}
	\begin{aligned}
		\dot{x}=&Ax+Bu+E\hat{C} \hat{v} \\
		y=&Cx\\
		e=&Cx+F\hat{C} \hat{v}
	\end{aligned}
\end{equation}
where $\hat{v}$ is known. To save  notation, without loss of generality, in this paper, we assume the exogenous signal $v$  in \eqref{lisys2} is known.
\end{rem}
}

	{\color{black}\begin{rem}
During the learning phase, the initial input takes the form $u^0=K_e \rho+\delta$, where $K_e$ is an initial exploring gain and $\delta$ is the exploration noise. As mentioned in \cite{Jiang2012, Bian2016} and  \cite{LH2024}, the fulfillment of the rank condition is  related to the persistent excitation property of the exploration noise $\delta$. A common choice of the exploration signal $\delta$ is a multitone sinusoidal function where the number of different frequencies depends on the number of unknown variables. It is worth mentioning that, by reducing the number of unknown variables,  it is easier to fulfill the rank condition required in our improved algorithm.
\end{rem}}

\subsection{Improved VI Algorithm}

{\color{black}It is possible to further reduce the computational complexity and the solvability condition of \eqref{linear equationR}  by the following procedure:

First, in \eqref{linear equationR} , letting $k=0$ gives
\begin{equation}\label{linear equationR0}
\begin{aligned}
\Psi^{\rho}\begin{bmatrix}
\text{vecs}(H_0^{\rho})\\
\text{vec}(E_\rho^TP_{0}^{\rho})
\end{bmatrix}=\Phi_{0}^{\rho}
\end{aligned}
\end{equation}

Under the rank condition \eqref{rank3}, solving \eqref{linear equationR0} gives $E_\rho$.
Then, the only unknown matrix in \eqref{linear equationR} is $H_k^{\rho}$.
%Since,
%for any symmetric matrix $H\in \mathbb{R}^{n\times n}$, there exists a constant matrix $M_n\in \mathbb{R}^{n^2\times \frac{n(n+1)}{2}}$ with full column rank such that $M_n\text{vecs}(H)=\text{vec}(H)$, we have
%\begin{align}\label{vinew111}
%		\textup{vec}(E^T_\rho P_k^{\rho})&=(I_{n_\rho}\otimes E^T_\rho )M_{n_\rho}\textup{vecs}(P_k^{\rho})
%	\end{align}

Thus, for $k\geq 1$, \eqref{linear equationR} can be simplified to the following:
\begin{align}\label{vinew22}
 I_{\rho \rho } \textup{vecs}(H_k^{\rho})={\Phi}^\rho_k - 2\Gamma_{\rho v}\text{vec}(E_\rho^TP_{k}^{\rho})
\end{align}

The matrix $I_{\rho \rho }$ has full column rank if
\begin{align}\label{rank4}
	\textup{rank}(I_{\rho \rho})=& \frac{n_{\rho}(n_{\rho}+1)}{2}
\end{align}		

Clearly, the rank condition \eqref{rank3} guarantees that the matrix  $I_{\rho \rho }$ has full column rank.

\begin{rem}
In our improved algorithm, for all $k \geq 1$, we have transformed the problem of solving linear equation \eqref{linear equationR} to that of solving \eqref{vinew22}.
Thus,  the number of  unknown variables is reduced by $q n_{\rho}$. Based on \eqref{linear equationR0} and \eqref{vinew22}, Algorithm 3 can be modified
as Algorithm 4.
\end{rem}

\begin{algorithm} \color{black}
\caption{The Improved Model-free VI Algorithm for Solving \eqref{AREMrho}} \label{vialg4}
\begin{algorithmic}[1]
       \State Apply any locally essentially bounded initial input $u^0$.
		Choose a proper $t_0 > 0$ such that $\|e_x\|$ is small enough.
		Collect data from $t_0$ until the rank condition \eqref{rank3} is satisfied.
		\State Choose $P_{0}^\rho \in \mathbb{R}^{n_\rho \times n_\rho}$ as an arbitrary positive definite matrix. Solve $H_0^{\rho}$ and $E_\rho$ from \eqref{linear equationR0}. $k,j\leftarrow 0$.
		\Loop
		\State Substitute $P_{k}^{\rho}$ into \eqref{vinew22} to solve $H_k^{\rho}$ and $K_{k}^{\rho}=-R^{-1}B_\rho^TP_{k}^{\rho}$.
		\State \begin{equation}\label{A4}
			\begin{aligned}
				\tilde{P}_{k+1}^{\rho} \gets P_{k}^{\rho}+\epsilon_k(H_k^{\rho}+Q_\rho-(K_{k}^{\rho})^TRK_{k}^{\rho})\\
			\end{aligned}
		\end{equation}
		\If{$\tilde{P}_{k+1}^{\rho} \notin \tilde{B}_j $}
		\State $P_{k+1}^{\rho} \gets P_{0}^{\rho}$. $j \gets j+1$.
		\ElsIf{$||\tilde{P}_{k+1}^{\rho}-P_{k}^{\rho}||/\epsilon_k<\varepsilon$}
		\State \textbf{return} $K_{k}^{\rho}$ as an approximation to $-R^{-1}B_\rho^TP_\rho^{*}$
		\Else
		\State $P_{k+1}^{\rho} \gets \tilde{P}_{k+1}^{\rho}$
		\EndIf
		\State $k \gets k+1$
		\EndLoop
	\end{algorithmic}
\end{algorithm}

}

\begin{rem}\label{rem6}
It is interesting to compare our control law (\ref{uz2}) with the one in \cite{Xie2023}, which takes the following form:
\begin{equation} \label{ref23}
	\begin{aligned}
	u & = K_x M_x \mbox{col} (\zeta,\zeta_v)+ K_v z\\
\dot{\zeta} & =(I_{m+p}\otimes \mathcal{A}) \zeta+ B_\zeta u+E_\zeta y \\
\dot{\zeta}_v&=(I_q \otimes \mathcal{A}) \zeta_v +(I_q \otimes b)v \\
\dot{z} & =G_1z+G_2e
\end{aligned}
\end{equation}
where $M_x = [M_1,\cdots ,M_{p+q+m}]$ where $M_i =  [D_0 f_i ,\cdots ,D_{n-1} f_i]$ with $f_i$ the $i^{th}$ column of $B$ for $i =1, \cdots, m$, $(i - m)^{th}$ column of $L$ for $i = m+1, \cdots, m+p$ and $(i - m-p)^{th}$ column of $E$ for $i = m+p+1, \cdots, m+p+q$.
{\color{black}First,  \eqref{ref23}  depends on $v$ and is  not in the output-feedback form. Second, \eqref{ref23}  contains an additional dynamics $\zeta_v$. Thus, the dimension of \eqref{ref23} is $n_\gamma=n(m+p+q)+n_z$, which is much higher than the dimension $n (p+m)+n_z$ of our control law.
Third, \eqref{ref23} is synthesized under the assumption that $(A, B)$ is controllable while we only assume that $(A, B)$ is stabilizable.
Fourth,  the algorithm of \cite{Xie2023} requires the fulfillment of the following rank condition:
\begin{align}\label{rankXie}
	\textup{rank}([I_{\gamma \gamma}, \Gamma_{\gamma u}, \Gamma_{\gamma v}])=& \frac{n_{\gamma}(n_{\gamma}+1)}{2}+(m+q)n_{\gamma}.
\end{align}
In contrast, we have implemented our control law with two much simplified algorithms, namely, Algorithms 3 and 4,  under much milder solvability condition \eqref{rank3}. }
%However, we claim that $\zeta_v$ converges to $Tv$ exponentially. At first, we consider the following Sylvester equation:
% \begin{align}\label{rankXie}
% TS=(I_q \otimes \mathcal{A})T+(I_q \otimes b)
% \end{align}
%admits a unique solution $T \in \mathbb{R}^{nq \times q}$ since the spectra of $S$ and $\mathcal{A}$ are disjoint.
%
%Then
%\begin{align}\label{rankXie}
%\dot{\zeta}_v-T\dot{v}=& \dot{\zeta}_v-TSv \notag \\
%=&(I_q \otimes \mathcal{A}) \zeta_v(t) +(I_q \otimes b)v(t)-TSv(t) \notag \\
%=&(I_q \otimes \mathcal{A})(\zeta_v-Tv)
%\end{align}
%
%Thus, $\zeta_v(t)=Tv(t)$ after $t \geq t_0$ for some finite $t_0$.
%
%Since $\xi=\mbox{col} (\zeta,\zeta_v,z)$, $[I_{\xi \xi}, \Gamma_{\xi u}, \Gamma_{\xi v}]=[\cdots, I_{\zeta_v \zeta_v}, \cdots, \Gamma_{\zeta_v v}, \cdots]$.
%$I_{\zeta_v \zeta_v}=I_{(Tv)(Tv)}=\Gamma_{(Tv)(Tv)}N_{nq}=\Gamma_{vv}(T\otimes T)^TN_{nq}$, where $N_{nq} \in \mathbb{R}^{(nq)^2 \times \frac{nq(nq+1)}{2}}$. $\Gamma_{\zeta_v v}=\Gamma_{(Tv) v}=\Gamma_{v v}(T \otimes I_q)^T$. Thus, the columns of $I_{\zeta_v \zeta_v}$ and $\Gamma_{\zeta_v v}$ are linearly dependent.
%
%Moreover, $\textup{rank}(I_{\zeta_v \zeta_v}) \leq \textup{rank}(\Gamma_{vv})\leq \frac{q(q+1)}{2}$,
%$\textup{rank}(\Gamma_{\zeta_v v}) \leq \textup{rank}(\Gamma_{vv})\leq \frac{q(q+1)}{2}$ which implies $I_{\zeta_v \zeta_v} \in \mathbb{R}^{s \times \frac{nq(nq+1)}{2}}$ and $\Gamma_{\zeta_v v} \in \mathbb{R}^{s \times nq^2}$ can not be full column rank when $n \geq 2$ or $q \geq 2$.
%In conclusion, $[I_{\xi \xi}, \Gamma_{\xi u}, \Gamma_{\xi v}]$ can not be full column rank.
\end{rem}

\begin{rem} \label{rem8}
A special case of the output regulation problem where $E=0$ was studied in \cite{Chen2022}. Even for this special case, the approach in this paper also offers two advantages. First, \cite{Chen2022} needs to estimate the whole state $\xi=\mbox{col} (x,z)$ of the augmented system  while we only estimate the state $x$ of the plant since  the state  $z$ of the internal model is known. Thus, the dimension of our control law is smaller than {\color{black}the control law in \cite{Chen2022} by $n_z(m+p)$.}
Second, since \cite{Chen2022} did not take advantage of the fact that  $B_\rho$ and hence  $K_k^{\rho}$ are known,   the approach of \cite{Chen2022} needs to solve a sequence of equations with larger dimension. In particular, the approach in \cite{Chen2022} requires the satisfaction of the  following rank condition:
{\color{black}\begin{align}\label{rankChen}
\textup{rank}([I_{\chi \chi}, \Gamma_{\chi u}, \Gamma_{\chi \nu}])=& \frac{n_{\chi}(n_{\chi}+1)}{2}+(m+p)n_\chi
\end{align}
where $\chi \in \mathbb{R}^{n_\chi}$ is the state of the augmented auxiliary  system proposed in \cite{Chen2022} with its dimension $n_\chi=(n+n_z)m+2(n+n_z)p$, and $\nu \in \mathbb{R}^p$ is the exploration noise.} The rank condition  in \eqref{rankChen} is more stringent than  \eqref{rank3}.  For a simple comparison, consider a case with $n=5, m=5, p=1, q=4, n_\rho=34, n_\chi=63$. {\color{black}Table \ref{tab2} compares the computational complexity of various algorithms proposed in \cite{Chen2022}, \cite{Xie2023} and this paper.}

\end{rem}

\begin{table}[!ht]\color{black}
\caption{Comparison of rank condition}\label{tab2}
\center
\begin{tabular}{|c|c|c|}\hline
Method & Equation labels & Unknown variables \\ \hline
Algorithm in \cite{Chen2022} & \eqref{rankChen} & 2394 \\
\hline
Algorithm in \cite{Xie2023} & \eqref{rankXie} & 1971 \\
\hline
Our Algorithm 3 & \eqref{rank3} &   731  \\
 \hline  		
Our Algorithm 4 & \eqref{rank4}&   595  \\
 \hline
 \end{tabular}
 \end{table}

%\begin{rem}
%	Our improved method can also be applied to \cite{liu2018adaptive} \cite{xie2022optimal} and \cite{gao2021reinforcement} to relax their solvability conditions and reduce their computing cost.
%\end{rem}

  %  Our improved VI-based data-driven algorithm is summarized as Algorithm \ref{vialg2},  whose convergence is guaranteed by Theorem \ref{converge}.

\section{LQR Problem for the Augmented auxiliary  System}\label{relaOP}

\subsection{The Relation between Problem \ref{StateLQR} and Problem \ref{OPFBLQR} }
It is known that, for any  ${Q}_{\xi} \geq 0$ such that $(Y, \sqrt{Q}_{\xi})$ is observable,  \eqref{riccati} admits the unique positive definite solution $P_{\xi}^{*}$ and
the  control law for solving Problem \ref{StateLQR} is
\begin{align}\label{lqrxi}
\bar{u}^*={K}_{\xi}^*\bar{\xi}
\end{align}
where ${K}_{\xi}^*=-R^{-1}J^T P_{\xi}^{*}$.
In what follows,  we will establish a relation between the solution ${P}_{\xi}^{*}$ to \eqref{riccati}  and the solution $P_\rho^{*} $ to \eqref{AREMrho}
with ${Q}_\rho=W^T{Q}_{\xi}W$, where $W=\mbox{blockdiag}(M, I_{n_z})$.
For this purpose,  let $\bar{C}= \begin{bmatrix} C & \bf{0} \\ \bf{0} & I_{n_z} \end{bmatrix}$ and $Q_{\xi} = \bar{C}^T \bar{Q}  \bar{C}$ where $\bar{Q} \in \mathbb{R}^{p+n_z} > 0$.
Then $\sqrt{Q_{\xi}} = \sqrt{\bar{Q}} \bar{C}$. Since Assumption \ref{ass2} implies the pair $(Y, \bar{C})$ is observable,
the pair $(Y, \sqrt{\bar{Q}} \bar{C} )$ is also observable.
We first show the following result.

\begin{lem}
Under Assumption \ref{ass2}, the pair $(A_\rho, \sqrt{\bar{Q}}\bar{C} W)$ is detectable.
\end{lem}

\begin{Prf}
%By definition, $A_\rho=\begin{bmatrix} A_\zeta & \bf{0}\\ G_2CM & G_1\end{bmatrix}$.
 Let $H \in \mathbb{R}^{n_z\times n_z}$ be any Hurwitz matrix,  and let $L_\rho=  \begin{bmatrix}-E_\zeta  & \bf{0}\\ -G_2  & H - G_1\end{bmatrix}  \sqrt{\bar{Q}}^{-1} $.
Then \begin{align*}
& A_\rho+L_\rho   \sqrt{\bar{Q}}\bar{C} W \\
=& \begin{bmatrix} I_{m+p}\otimes \mathcal{A} + E_\zeta CM & \bf{0}\\
G_2 C M  & G_1 \end{bmatrix} + \begin{bmatrix}-E_\zeta CM & \bf{0}\\ -G_2 CM & (H-G_1) \end{bmatrix} \\
=&\begin{bmatrix} I_{m+p}\otimes \mathcal{A} & \bf{0}\\
\bf{0} & H\end{bmatrix}
\end{align*}

Since $\mathcal{A}$ and $H$ are Hurwitz, $A_\rho+L_\rho   \sqrt{\bar{Q}}\bar{C} W$ is Hurwitz. Thus, the pair $(A_\rho, \sqrt{\bar{Q}}\bar{C} W)$ is detectable.
\end{Prf}

We are ready to establish the main result of this subsection.
\begin{thm} \label{thm4}
Under Assumptions \ref{ass1}-\ref{ass5}, for any  ${\bar{Q}} > 0$ and $Q_{\xi}= \bar{C}^T\bar{Q} \bar{C}$, the algebraic Riccati equation \eqref{AREMrho}
 with ${Q}_\rho=W^T {Q}_{\xi}W$ admits a unique positive semidefinite solution ${P}_\rho^{*}=W^T {P}_{\xi}^{*}W$, where ${P}_{\xi}^{*}$ is the unique positive definite solution to \eqref{riccati}. Thus, the solution to Problem \ref{OPFBLQR} is given by
 \begin{align}\label{lqrrho}
\bar{u}^*=K_\rho^*\bar{\rho}
\end{align}
where $K_\rho^*=-R^{-1}B_\rho^T {P}_\rho^{*} = -R^{-1}B_\rho^T W^T {P}_{\xi}^{*}W = {K}_{\xi}^* W$.
%As a result, the control law that solves Problem \ref{OPFBLQR}  asymptotically approaches the state-feedback control law that solves Problem \ref{StateLQR}.
\end{thm}

\begin{Prf}
Since $(A_\rho, B_\rho)$ is stabilizable and $(A_\rho, \sqrt{\bar{Q}}\bar{C} {W})$ is detectable,
the algebraic Riccati equation \eqref{AREMrho}
 with ${Q}_\rho=W^T{Q}_{\xi}W$ admits a unique positive semidefinite solution ${P}_\rho^{*}$.
 We only  need to show that ${P}_\rho^{*}=W^T {P}_{\xi}^{*}W$, where ${P}_{\xi}^{*}$ is the unique positive definite solution to \eqref{riccati}.
 Using  \eqref{rela1}-\eqref{rela3} shows that
\begin{equation}\label{WA}
\begin{aligned}
WA_\rho=&\begin{bmatrix}
M &\bf{0} \\
\bf{0} & I_{n_z}
\end{bmatrix}
\begin{bmatrix}
A_\zeta & \bf{0}\\
G_2CM & G_1
\end{bmatrix}
=\begin{bmatrix}
MA_\zeta &\bf{0} \\
G_2CM & G_1\\
\end{bmatrix}\\
=&\begin{bmatrix}
A &\bf{0} \\
G_2C & G_1\\
\end{bmatrix}
\begin{bmatrix}
M &\bf{0} \\
\bf{0} & I_{n_z}
\end{bmatrix}= Y W
\end{aligned}
\end{equation}

\begin{equation}\label{WB}
\begin{aligned}
WB_\rho=&\begin{bmatrix}
M &\bf{0} \\
\bf{0} & I_{n_z}
\end{bmatrix}
\begin{bmatrix}
B_\zeta \\
\bf{0}
\end{bmatrix}
=\begin{bmatrix}
MB_\zeta \\
\bf{0}\\
\end{bmatrix}
=\begin{bmatrix}
B \\
\bf{0}\\
\end{bmatrix}=J
\end{aligned}
\end{equation}

Using \eqref{WA} and \eqref{WB}, it can be verified that

\begin{equation}
\begin{aligned}
&A_\rho^T(W^T {P}_{\xi}^{*}W)+(W^T {P}_{\xi}^{*}W)A_\rho+{Q}_\rho\\
&-(W^T {P}_{\xi}^{*}W)B_\rho R^{-1}B_\rho^T(W^T {P}_{\xi}^{*}W)\\
=&W^T {Y}^T {P}_{\xi}^{*}W+W^T {P}_{\xi}^{*} Y W+W^T {Q}_{\xi}W\\
&-W^T {P}_{\xi}^{*}J R^{-1} J ^T {P}_{\xi}^{*}W\\
=&W^T(Y^T{P}_{\xi}^{*}+{P}_{\xi}^{*}Y+{Q}_{\xi}-{P}_{\xi}^{*} J R^{-1} J^T {P}_{\xi}^{*})W\\
=&0
\end{aligned}
\end{equation}

Thus, $W^T {P}_{\xi}^{*}W$ is a positive semidefinite solution to the algebraic Riccati equation \eqref{AREMrho}. Since ${P}_\rho^{*}$ is unique, ${P}_\rho^{*}=W^T{P}_{\xi}^{*}W$.

Then substituting ${P}_\rho^{*}=W^T{P}_{\xi}^{*}W$ into ${K}_\rho^* $ and using \eqref{WB} gives
\begin{equation}
\begin{aligned}
{K}_\rho^*=&-R^{-1}B_\rho^T {P}_\rho^{*}\\
=&-R^{-1}B_\rho^T W^T {P}_{\xi}^{*}W\\
=&-R^{-1} J^T {P}_{\xi}^{*}W\\
=&{K}_{\xi}^*W
\end{aligned}
\end{equation}

\end{Prf}

\subsection{The VI Method for Problem \ref{OPFBLQR}}
{\color{black}
Theorem \ref{thm4} has established the relation for the control law \eqref{lqrxi} that solves Problem 2 and the control law \eqref{lqrrho} that solves Problem 4, that is,  if ${Q}_\rho=W^T {Q}_{\xi}W$, then   ${P}_\rho^{*}=W^T {P}_{\xi}^{*}W$ and $K_\rho^*= {K}_{\xi}^* W$. However,
this particular ${Q}_\rho=W^T {Q}_{\xi}W$ is not feasible because $M$ and hence $W=\mbox{blockdiag}(M, I_{n_z})$ are unknown. Nevertheless, for the special case
where $E=\bf{0}$, it is possible to circumvent this difficulty. For this purpose, we first establish the following result.

\begin{lem} \label{lem14}
Under Assumptions \ref{ass1}-\ref{ass5}, if $E = \bf{0}$, then,  with $\bar{Q} = \begin{bmatrix}
Q_y & \bf{0} \\
\bf{0} & Q_z
\end{bmatrix}$ where $Q_y >0$ and $Q_z >0$,
$|{\rho}|_{{Q}_\rho}$ approaches $( |y|_{Q_y} + |z|_{Q_z})$ exponentially.
\end{lem}

\begin{Prf}

First note that, if $E = \bf{0}$, then $X'=\bf{0}$ in \eqref{SylX}. By \eqref{x4}, $M\zeta(t)$ approaches $x(t)$ exponentially.
Next, with $\bar{Q} = \begin{bmatrix}
Q_y & \bf{0} \\
\bf{0} & Q_z
\end{bmatrix}$,
${Q}_{\xi} = \begin{bmatrix}
C^T Q_y C & \bf{0} \\
\bf{0} & Q_z
\end{bmatrix}$. Since  ${\xi} = \col ({x}, {z})$,
$|{\xi}|_{Q_{\xi}} = |C {x}| _{Q_y}+ |{z}|_{Q_z} = |{y}| _{Q_y}+ |{z}|_{Q_z}$ where ${y} = C {x}$.
On the other hand, since ${Q}_{\rho} = W^T Q_\xi W$,
 $|{\rho}|_{{Q}_\rho}= (M \zeta)^T C^T Q_y C M \zeta  + |z|_{Q_z}$. Since  $M\zeta(t)$ approaches $x(t)$ exponentially,  $|{\rho}|_{{Q}_\rho}$ approaches $( |y|_{Q_y} + |z|_{Q_z})$ exponentially.
\end{Prf}

We now propose a model-free VI method to calculate ${P}_\rho^{*}$ and $K_\rho^*$.
For this purpose,  note that if $E = \bf{0}$,
 \eqref{rhoz5} reduces to the following
\begin{align}\label{rhoz6}
\dot{\rho}=A_\rho \rho+B_\rho u+\tilde{E}_\rho v.
\end{align}
 where $\tilde{E}_\rho=\begin{bmatrix}\bf{0} \\ G_2F \end{bmatrix}$.

}

Integrating $\frac{ d (\rho^T (t) {{P}^\rho_k} \rho (t))}{dt}$ along the solution of  \eqref{rhoz6}
gives
\begin{equation}\label{RZIRLtil2}
\begin{aligned}
&|\rho(t+\delta t)|_{P_{k}^{\rho}}-|\rho(t)|_{P_{k}^{\rho}}\\
=&\int_{t}^{t+\delta t} \left( |\rho|_{H_k^{\rho}}  + 2u^TB_\rho^TP_{k}^{\rho}\rho +2v^T\tilde{E}_\rho^TP_{k}^{\rho} \rho \right) d\tau
\end{aligned}
\end{equation}
\eqref{RZIRLtil2} can be rewritten to the following:
\begin{equation}\label{RZIRLtil}
\begin{aligned}
&|\rho(t+\delta t)|_{P_{k}^{{\rho}}}-|\rho(t)|_{P_{k}^{{\rho}}}+\int_{t}^{t+\delta t} \left( |\rho|_{{Q}_\rho} \right) d\tau\\
=&\int_{t}^{t+\delta t} \left( |\rho|_{\bar{H}_k^{\rho}}+2u^TB_\rho^TP_{k}^{{\rho}} \rho +2v^T\tilde{E}_\rho^TP_{k}^{\rho} \rho \right) d\tau\\
\end{aligned}
\end{equation}
where $\bar{H}_k^{\rho}= {H}_k^{\rho}+{Q}_\rho$.

{\color{black} By Lemma \ref{lem14},  with  ${Q}_{\xi} = \begin{bmatrix}
C^T Q_y C & \bf{0} \\
\bf{0} & Q_z
\end{bmatrix}$, $|\rho|_{{Q}_\rho}$ approaches $( |y|_{Q_y} + |z|_{Q_z})$ exponentially.   Thus, } \eqref{RZIRLtil} can be approximated by the following equation:
\begin{equation}\label{RZIRLtil3}
\begin{aligned}
&|\rho(t+\delta t)|_{P_{k}^{{\rho}}}-|\rho(t)|_{P_{k}^{{\rho}}} +\int_{t}^{t+\delta t} \left( |y|_{Q_y} + |z|_{Q_z} \right) d\tau\\
&=\int_{t}^{t+\delta t} \left( |\rho|_{\bar{H}_k^{\rho}} + 2u^TB_\rho^T P_{k}^{{\rho}} \rho +2v^T\tilde{E}_\rho^TP_{k}^{\rho} \rho \right) d\tau
\end{aligned}
\end{equation}

\eqref{RZIRLtil3} implies
\begin{equation}\label{linear equationRtil}
\begin{aligned}
\tilde{\Psi}^{\rho}\begin{bmatrix}
\text{vecs}(\bar{H}_k^{\rho})\\
\text{vec}(\tilde{E}_\rho^TP_{k}^{\rho})
\end{bmatrix}=\tilde{\Phi}_{k}^{\rho}
\end{aligned}
\end{equation}
where
\begin{equation}
\begin{aligned}
\tilde{\Psi}^{\rho}=&[I_{\rho \rho}, 2\Gamma_{\rho v}]\\
\tilde{\Phi}_{k}^{\rho}=&\delta_{\rho}\text{vecs}(P_{k}^{\rho})+I_{yy}\text{vecs}(Q_y)+I_{zz}\text{vecs}(Q_z)\\
&-2\Gamma_{\rho (B_\rho u)}\text{vec}(P_{k}^{\rho})
\end{aligned}
\end{equation}

The solvability of \eqref{linear equationRtil} is guaranteed by the following lemma.
\begin{lem} \label{rankVI3til}
The matrix $\tilde{\Psi}^{\rho}$ has full column rank if
\begin{align}\label{rank3til}
\textup{rank}([I_{\rho \rho}, \Gamma_{\rho v}])=& \frac{n_{\rho}(n_{\rho}+1)}{2}+qn_{\rho}
\end{align}
\end{lem}

Similar to Section IV.C, we can further simplify \eqref{linear equationRtil} as follows:

{\color{black}
First, in \eqref{linear equationRtil}, letting $k=0$ gives
\begin{equation}\label{linear equationRtil0}
	\begin{aligned}
		\tilde{\Psi}^{\rho}\begin{bmatrix}
			\text{vecs}(\bar{H}_0^{\rho})\\
			\text{vec}(\tilde{E}_\rho^TP_{0}^{\rho})
		\end{bmatrix}=\tilde{\Phi}_{0}^{\rho}
	\end{aligned}
\end{equation}

Under the rank condition \eqref{rank3til}, solving \eqref{linear equationRtil0} gives $\tilde{E}_\rho$.
Then, the only unknown matrix in \eqref{linear equationRtil} is $\bar{H}_k^{\rho}$.

Thus, for $k\geq 1$, \eqref{linear equationRtil} can be simplified to the following:
\begin{align}\label{vinew22til}
I_{\rho \rho } \textup{vecs}(\bar{H}_k^{\rho})=\tilde{\Phi}_{k}^{\rho}-2\Gamma_{\rho v} \text{vec}(\tilde{E}_\rho^TP_{k}^{\rho})
\end{align}

Clearly, the rank condition \eqref{rank3til} guarantees that the matrix  $I_{\rho \rho }$ has full column rank. Based on \eqref{linear equationRtil0} and \eqref{vinew22til}, Algorithm 5 can be modified as Algorithm 6. }

%Using the technique given in \cite{LH2024}, the rank condition in \eqref{rank3til} is weaker than \eqref{rank3}. For a simple comparison, consider a case with $n=5, m=5, p=1, q=4, n_\rho=34$. Table \ref{tab1} compares the computational complexity of the original VI Algorithm, and the improved algorithm.

%\begin{table}[!ht]
%\caption{Comparison of rank condition}\label{tab1}
%\center
%\begin{tabular}{|c|c|c|}\hline
%Method & Equation & Unknown variables \\ \hline
%Original algorithm & \eqref{rank3} & 731 \\
%\hline
 %Improved algorithm & \eqref{rank3til} & 595   \\
% \hline  		
% \end{tabular}
% \end{table}

\begin{algorithm}
\caption{The  VI Algorithm for Solving \eqref{AREMrho} with ${Q}_\rho=W^T {Q}_{\xi}W$ and $E=0$} \label{vialg5}
\begin{algorithmic}[1]
\State Apply any locally essentially bounded initial input $u^0$.
Choose a proper $t_0 > 0$ such that $\|e_x\|$ is small enough.
Collect data from $t_0$ until the rank condition \eqref{rank3til} is satisfied.
\State Choose $P_{0}^{{\rho}}=(P_{0}^{{\rho}})^T \geq 0$. $k,j \gets 0$
\Loop
\State Substitute $P_{k}^{{\rho}}$ into \eqref{linear equationRtil} to solve $\bar{H}_k^{{\rho}}$ and $K_{k}^{{\rho}}=-R^{-1}B_\rho^TP_{k}^{\rho}$.
\State \begin{equation}\label{A5}
\begin{aligned}
\tilde{P}_{k+1}^{{\rho}} \gets P_{k}^{{\rho}}+\epsilon_k(\bar{H}_k^{{\rho}}-(K_{k}^{{\rho}})^TRK_{k}^{{\rho}})\\
\end{aligned}
\end{equation}
\If{$\tilde{P}_{k+1}^{{\rho}} \notin \tilde{B}_j $}
    \State $P_{k+1}^{{\rho}} \gets P_{0}^{{\rho}}$. $j \gets j+1$.
\ElsIf{$||\tilde{P}_{k+1}^{{\rho}}-P_{k}^{{\rho}}||/\epsilon_k<\varepsilon$}
    \State \textbf{return} $K_{k}^{{\rho}}$ as an approximation to $-R^{-1}B_\rho^TP_\rho^{*}$
\Else
	\State $P_{k+1}^{{\rho}} \gets \tilde{P}_{k+1}^{{\rho}}$
\EndIf
\State $k \gets k+1$
\EndLoop
\end{algorithmic}
\end{algorithm}

\begin{algorithm} 
	\caption{The Improved VI Algorithm for Solving \eqref{AREMrho} with ${Q}_\rho=W^T {Q}_{\xi}W$ and $E=0$} \label{vialg6}
	\begin{algorithmic}[1]
		\State Apply any locally essentially bounded initial input $u^0$.
		Choose a proper $t_0 > 0$ such that $\|e_x\|$ is small enough.
		Collect data from $t_0$ until the rank condition \eqref{rank3til} is satisfied.
	\State Choose $P_{0}^\rho \in \mathbb{R}^{n_\rho \times n_\rho}$ as an arbitrary positive definite matrix. Solve $\bar{H}_0^{\rho}$ and $\tilde{E}_\rho$ from \eqref{linear equationRtil0}. $k,j\leftarrow 0$.
		\Loop
		\State Substitute $P_{k}^{{\rho}}$ into \eqref{vinew22til} to solve $\bar{H}_k^{{\rho}}$ and $K_{k}^{{\rho}}=-R^{-1}B_\rho^TP_{k}^{\rho}$.
		\State \begin{equation}\label{A6}
			\begin{aligned}
				\tilde{P}_{k+1}^{{\rho}} \gets P_{k}^{{\rho}}+\epsilon_k(\bar{H}_k^{{\rho}}-(K_{k}^{{\rho}})^TRK_{k}^{{\rho}})\\
			\end{aligned}
		\end{equation}
		\If{$\tilde{P}_{k+1}^{{\rho}} \notin \tilde{B}_j $}
		\State $P_{k+1}^{{\rho}} \gets P_{0}^{{\rho}}$. $j \gets j+1$.
		\ElsIf{$||\tilde{P}_{k+1}^{{\rho}}-P_{k}^{{\rho}}||/\epsilon_k<\varepsilon$}
		\State \textbf{return} $K_{k}^{{\rho}}$ as an approximation to $-R^{-1}B_\rho^TP_\rho^{*}$
		\Else
		\State $P_{k+1}^{{\rho}} \gets \tilde{P}_{k+1}^{{\rho}}$
		\EndIf
		\State $k \gets k+1$
		\EndLoop
	\end{algorithmic}
\end{algorithm}

\section{Numerical Examples}\label{eg}
Consider the linear system \eqref{lisys2} with
\begin{equation}\nonumber
\begin{aligned}
&{\color{black}A=\begin{bmatrix}
0 & 1  & 0 \\ 
0 & 0  & 0 \\ 
0 & 0  &-1
\end{bmatrix},
B=\begin{bmatrix}
0\\
1\\
0
\end{bmatrix},}
S=\begin{bmatrix}
0 & 1\\
-1 & 0
\end{bmatrix},\\
&{\color{black}C=\begin{bmatrix}
1 & 2 & 3
\end{bmatrix},
F=\begin{bmatrix}
0.5 & -0.8
\end{bmatrix}}.
\end{aligned}
\end{equation}
{\color{black}Since $(A, B)$ is stabilizable but not controllable, the existing methods in \cite{Chen2022} and \cite{Xie2022} cannot handle this example. 
The internal model for  this example is as follows: }
\begin{equation}\nonumber
\begin{aligned}
G_1=\begin{bmatrix}
0 & 1\\
-1 & 0
\end{bmatrix},
G_2=\begin{bmatrix}
0 \\
1
\end{bmatrix}.
\end{aligned}
\end{equation}
 {\color{black}Let $
	L=\begin{bmatrix}
		-523 & 210 & 40
	\end{bmatrix}^T$, which is such that $
	\lambda(A-LC)=\{-5,-6,-7\}.$
Also, we have 
\begin{equation}
\begin{aligned}
		M=&\begin{bmatrix}
M_u & M_y
\end{bmatrix}
\end{aligned}
\end{equation}
where
\begin{equation}
\begin{aligned}
M_u=&10^3 \times \begin{bmatrix}
1.1670  &  1.0470  &     0 \\
-0.5230 &  -0.4020 &   0.0010 \\
-0.0400 &  -0.0800  &       0
\end{bmatrix}\\
M_y=&10^3 \times \begin{bmatrix}
 0.2100 &  -0.3130 &  -0.5230\\
0.0000 &   0.2100 &   0.2100\\
 0.0000 &  -0.0000 &   0.0400
\end{bmatrix}
\end{aligned}
\end{equation}}

\subsection{The scenario where $E=\bf{0}$}
{\color{black}
This scenario will be handled by Algorithm \ref{vialg6}.  The initial conditions are randomly chosen as $v(0)=\mbox{col}(1, 0.8)$, $x(0)=\mbox{col}(1, 2, -0.8)$, $z(0)=\mbox{col}(0, 0)$, $\zeta(0)=\mbox{col}(0, 0, 0, 0, 0, 0)$. The simulation is performed with the exploration initial input $u^0=K_0 \zeta +\delta_1$, where $K_0=\begin{bmatrix}
10 & 8 & 0 & 0 & -4 & -4
\end{bmatrix}$ and  $\delta_1=10(\sin(4t)+\sin(10t)+\sin(9t)-\sin(2t)-\sin(6t))$ is an exploration noise. The observers take about four seconds to converge. Let $t_j=4+0.2j$ with $j=0,1,\cdots,s$ and $s=120$.  Collecting the data  from $4\sec$ to $28\sec$  guarantees \eqref{rank3til}.
For implementing  Algorithm \ref{vialg6}, various parameters are set as $P_{0}^{\rho}=0.01I_8$, $B_{\rho}=\mbox{col} (0, 0, 1, 0, 0, 0, 0, 0)$, $\tilde{B}_j=\{P \in \mathcal{P}^8: |P|\leq 1000(j+20)\}$, $\epsilon_k=8/(k+10)$, $\varepsilon=0.05$, ${Q}_y =1$, ${Q}_z = I_{2}$ and $R=1$.
Then, it takes $8771$ iterations for the convergence criterion $||\tilde{P}_{k+1}^{{\rho}}-P_{k}^{{\rho}}||/\epsilon_k<\varepsilon$ to be satisfied. The comparison of $P_{k}^{{\rho}}$ and $K_{k}^{{\rho}}$ with the desired values ${P}_\rho^{*}$ and ${K}_\rho^{*}$ are shown in Fig. \ref{E=0}.

%\begin{figure}
%	\centering
%     \includegraphics[width=1.06\linewidth]{tilP}
%	\caption{Comparison of $P_{k}^{{\rho}}$, $K_{k}^{{\rho}}$ and the actual values ${P}_\rho^{*}$, ${K}_\rho^{*}$}
%	\label{tildeP}
%\end{figure}

\begin{figure}
	\centering
   \hspace{-0.6cm}
	\subfigure
		\centering
		\includegraphics[width=0.54\linewidth]{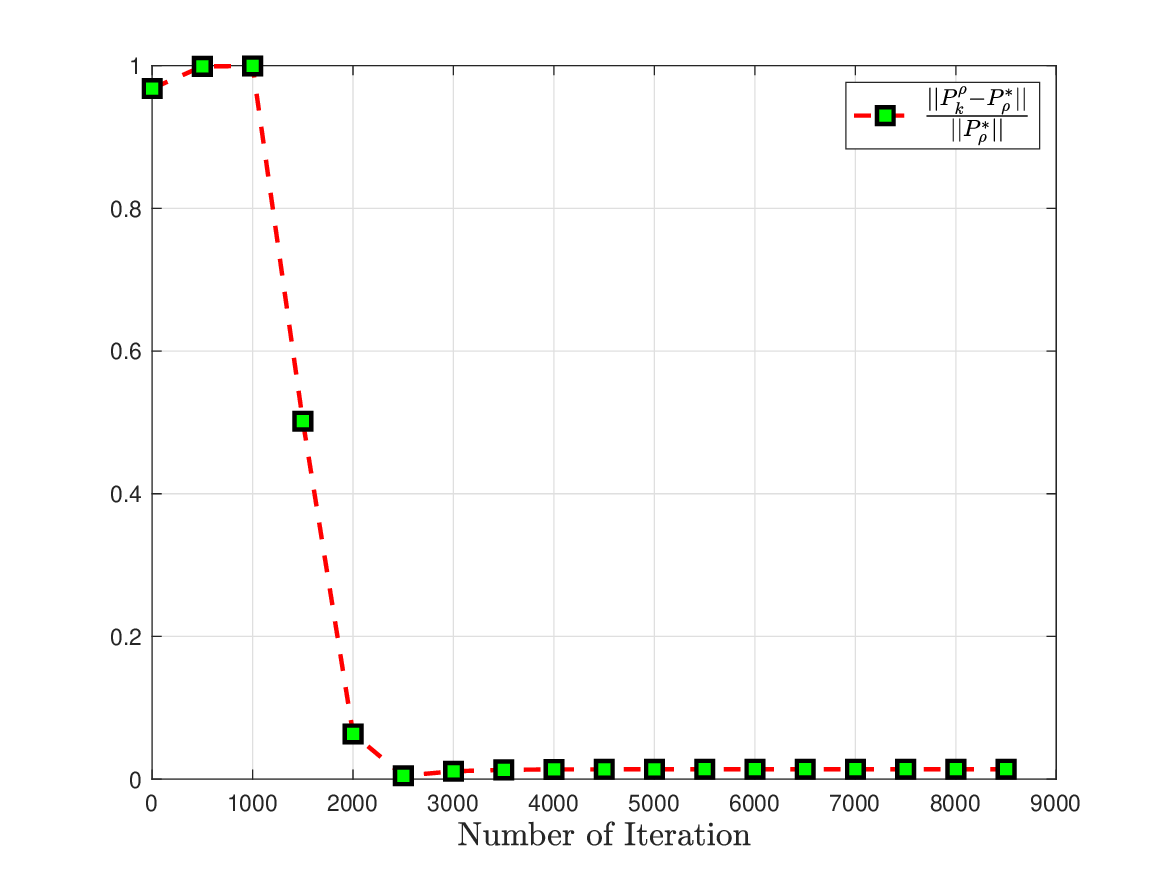}
		\label{p1}
	\hspace{-0.5cm}
	\subfigure
		\centering
		\includegraphics[width=0.54\linewidth]{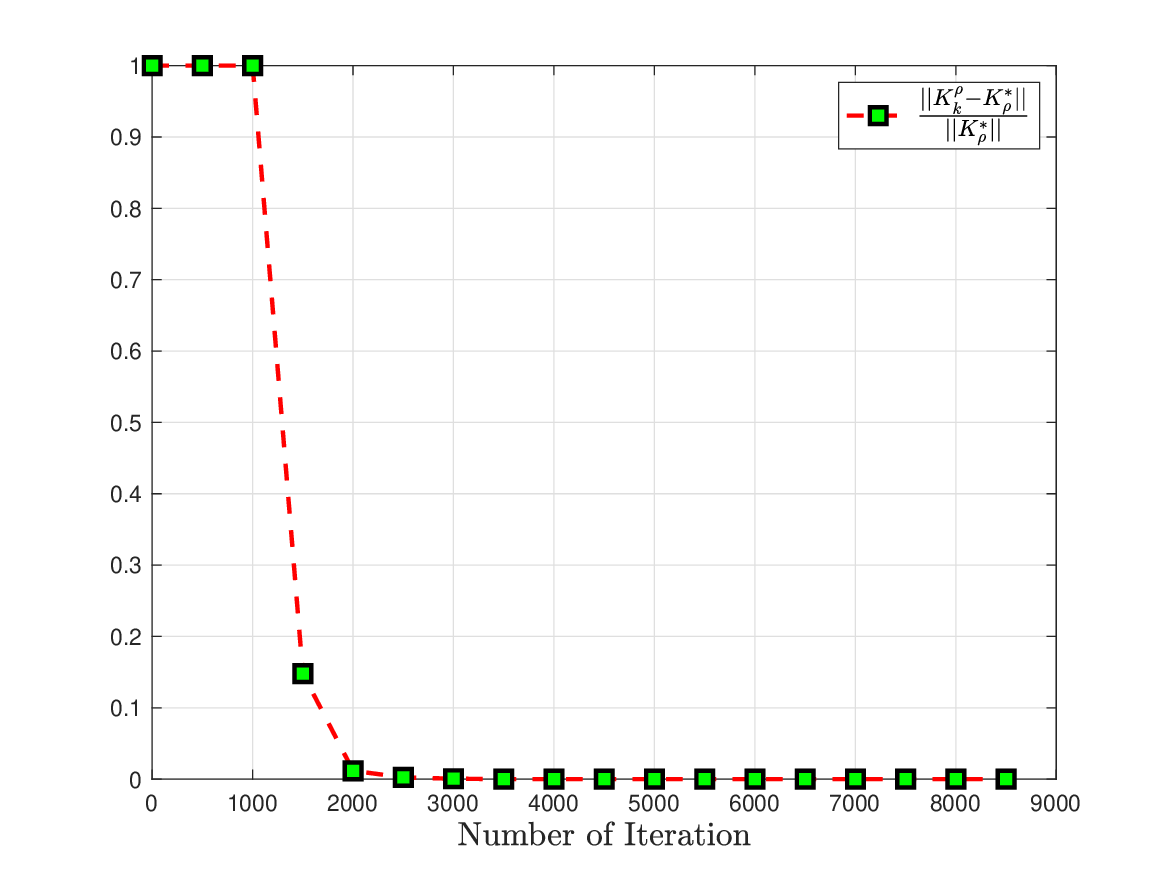}
		\label{k1}
	\caption{Comparison of $P_{k}^{{\rho}}$, $K_{k}^{{\rho}}$ with the desired values ${P}_\rho^{*}$, ${K}_\rho^{*}$}
	\label{E=0}
\end{figure}

Control law $u=K_{8771}^{{\rho}} \rho$ starts to operate at $t=28\sec$. Fig. \ref{e1} shows the error between the measurement output $y$ and the reference signal $-Fv$ after applying the control law. As expected, the controller works satisfactorily.
\begin{figure}
	\centering
      \includegraphics[width=1.0\linewidth]{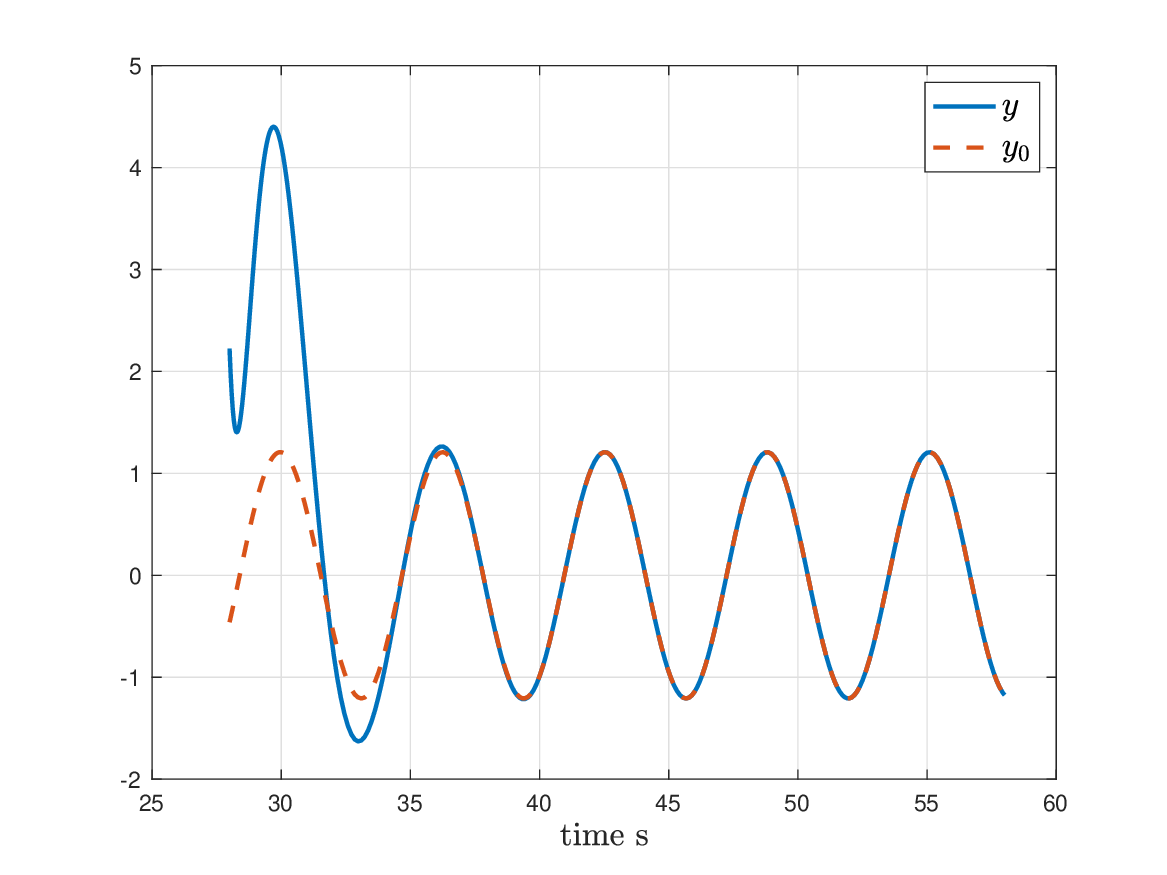}
	\caption{Plots of the system output $y$ and the reference signals $-Fv$ after control policy updated.}
	\label{e1}
\end{figure}}

\subsection{The scenario where $E\neq \bf{0}$}
{\color{black} We apply Algorithm \ref{vialg4} to this scenario. 
Let
$E=\begin{bmatrix}
2 & 0\\
0 & 1\\
3 & 6
\end{bmatrix}$.

The initial conditions are randomly chosen as $v(0)=\mbox{col}(1, 0.8)$, $x(0)=\mbox{col}(1, 2, -0.8)$, $z(0)=\mbox{col}(0, 0)$, $\zeta(0)=\mbox{col}(0, 0, 0, 0, 0, 0)$. The simulation is performed with the exploration initial input $u^0=K_0\zeta +\delta$, where $K_0=\begin{bmatrix}
10 & 8 & 0 & 0 & -4 & -4
\end{bmatrix}$ and $\delta=10(\sin(4t)+\sin(9t)+\sin(10t)-\sin(2t)-\sin(6t))$ is an exploration noise. The observers takes about four seconds to  converge. Let $t_j=4+0.2j$ with $j=0,1,\cdots,s$ and $s=120$. Collecting the data  from $4\sec$ to $28\sec$ guarantees \eqref{rank3}.

To implement Algorithm \ref{vialg4}, we set various parameters to be $P_{0}^{\rho}=0.1I_8$, $B_{\rho}=\mbox{col} (0, 0, 1, 0, 0, 0, 0, 0)$, $\tilde{B}_j=\{P \in \mathcal{P}^8: |P|\leq 1000(j+20)\}$, $\epsilon_k=20/(k+4000)$, $\varepsilon=0.01$, $Q_\rho=I_8$ and $R=1$. Then, it takes $10602$ iterations for the convergence criterion $||\tilde{P}_{k+1}^{\rho}-P_{k}^{\rho}||/\epsilon_k<\varepsilon$ to be satisfied. The comparison of $P_{k}^{{\rho}}$ and $K_{k}^{{\rho}}$ with the desired values ${P}_\rho^{*}$ and ${K}_\rho^{*}$ are shown in Fig. \ref{Eneq0}.

\begin{figure}
	\centering
	\hspace{-0.6cm}
	\subfigure
	\centering
	\includegraphics[width=0.54\linewidth]{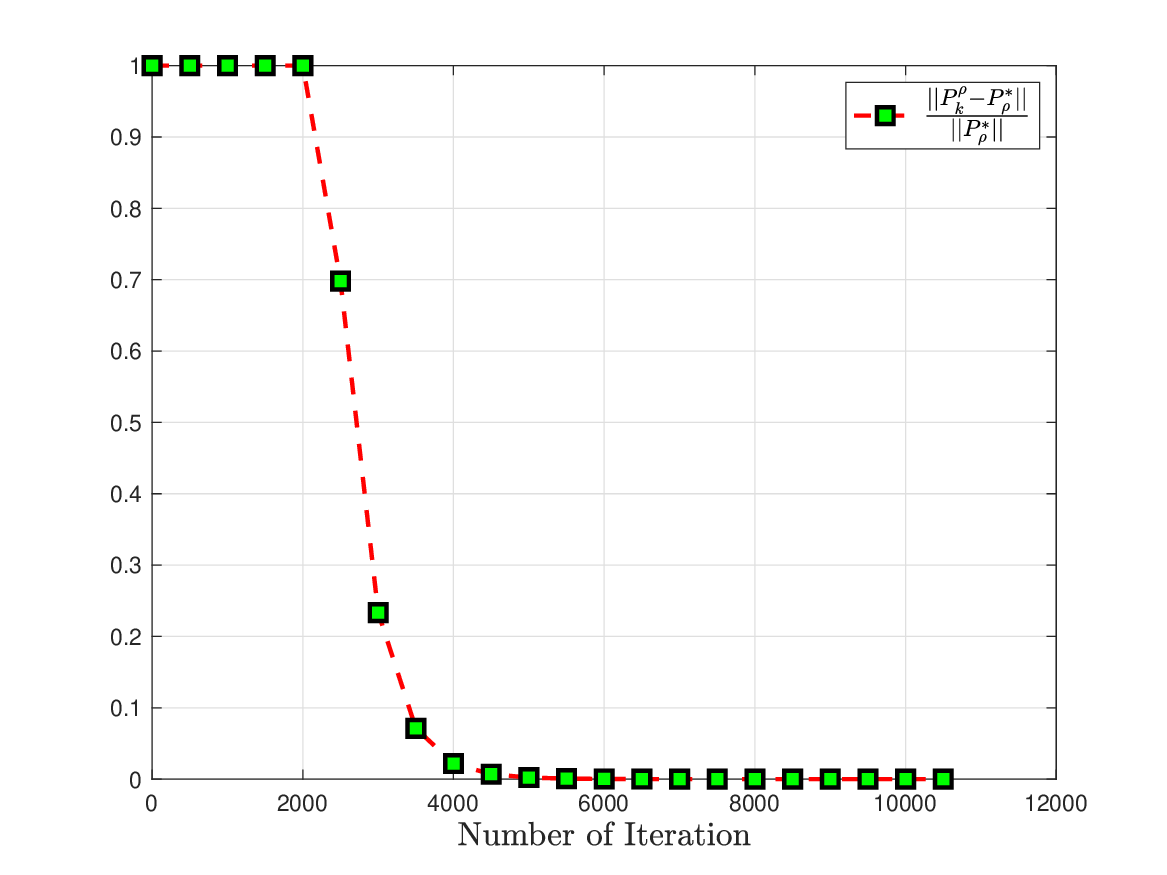}
	\label{p2}
	\hspace{-0.5cm}
	\subfigure
	\centering
	\includegraphics[width=0.54\linewidth]{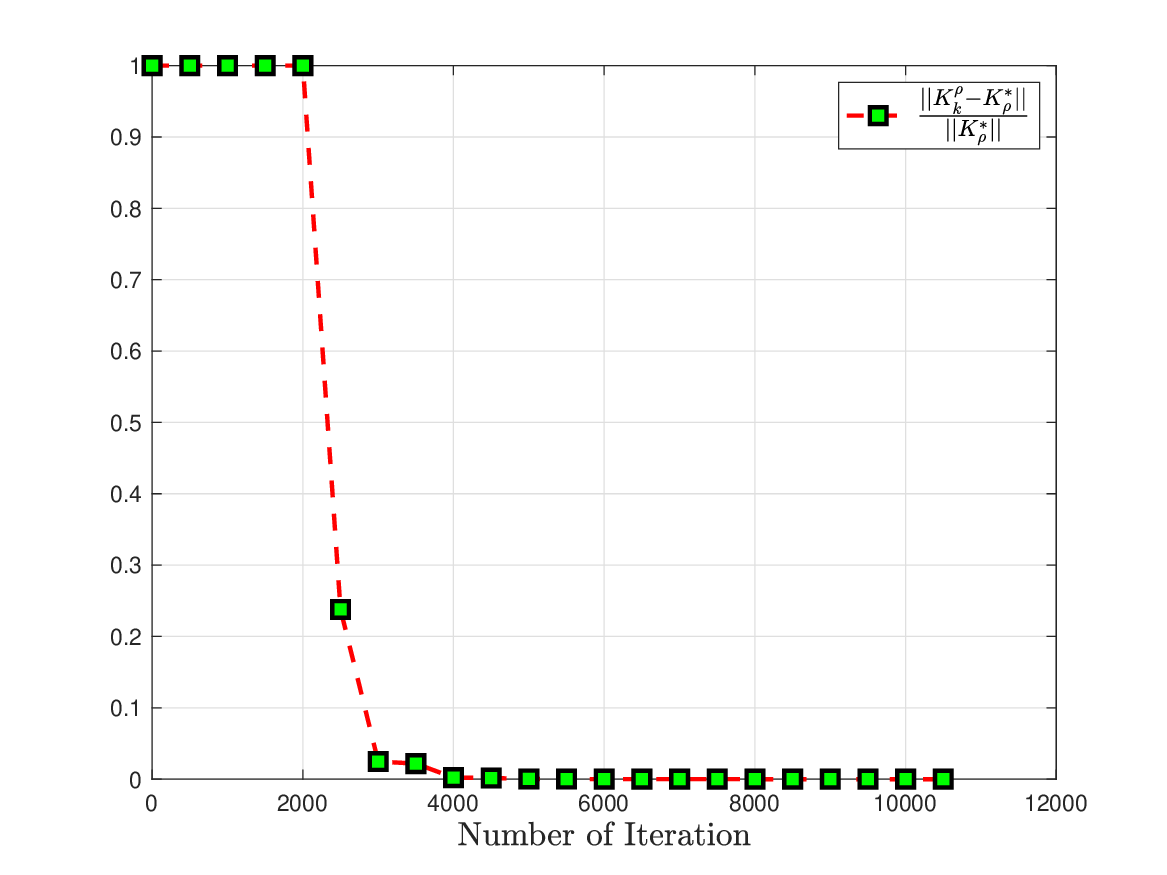}
	\label{k2}
	\caption{Comparison of $P_{k}^{{\rho}}$, $K_{k}^{{\rho}}$ with the desired values ${P}_\rho^{*}$, ${K}_\rho^{*}$}
	\label{Eneq0}
\end{figure}

Control law $u=K_{10602}^{\rho} \rho$ starts to operate at $t=28\sec$. Fig. \ref{e2} shows error between the measurement output $y$ and the reference signal $-Fv$ after applying the control law. As expected, the controller works satisfactorily.
\begin{figure}
	\centering
      \includegraphics[width=1.0\linewidth]{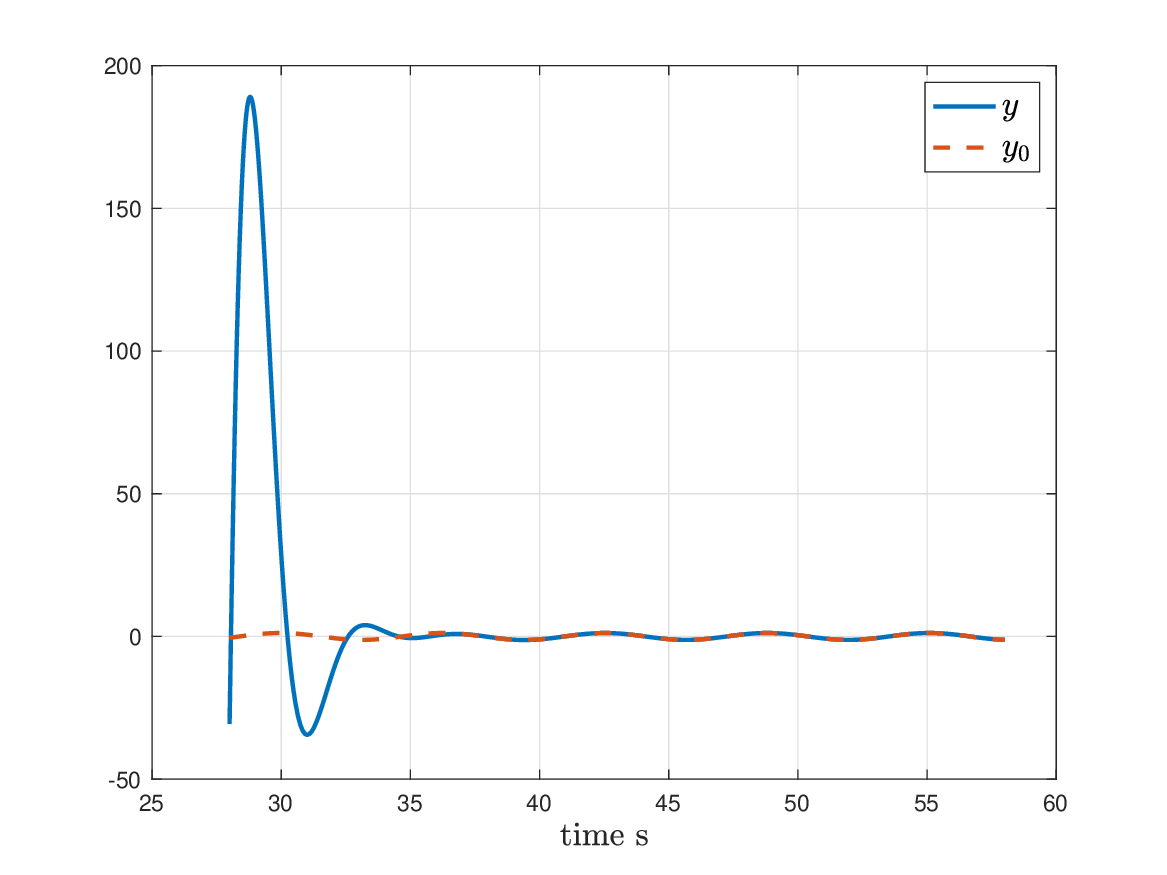}
	\caption{Plots of the system output $y$ and the reference signals $-Fv$ after control policy updated.}
	\label{e2}
\end{figure}}

\begin{rem}
{\color{black}	As pointed out in Remark \ref{rem6},   the algorithm of \cite{Xie2023} requires the fulfillment of the rank condition \eqref{rankXie}. 
	For this example, $n_\gamma=n(m+p+q)+n_z=14$. Thus, the rank condition \eqref{rankXie} requires 
	\begin{align}\label{rankXie2}
		\textup{rank}([I_{\gamma \gamma}, I_{\gamma u}, I_{\gamma v}])=147
	\end{align}
However, even if we let $s = 200$}, 	$\textup{rank}([I_{\gamma \gamma}, I_{\gamma u}, I_{\gamma v}])= 63 < 147$. 
\end{rem}

\section{Conclusions}\label{conclusion}
In this paper, we have studied the output regulation problem for unknown linear systems by the  data-driven output-based approach via value iteration. We first developed a novel output-feedback control law that does not explicitly rely on the observer gain to solve the output regulation problem.
Then we have further shown that the data-driven approach for designing an output-feedback control law for the given plant can be converted to the data-driven approach for designing a state-feedback control law for the augmented auxiliary  system and {\color{black}have thus obtained} a systematic data-driven approach for solving the output regulation problem for the unknown linear system via value iteration.
Finally, we have established a relation between the data-driven state-feedback control law and the data-driven output-feedback control law  in the LQR sense.
Compared with the existing results, our control law only depends on the output of the plant and {\color{black}deals with both the asymptotic tracking and disturbance rejection.}
Also, our algorithm has significantly reduced the computational cost and weakened the solvability condition.

\ifCLASSOPTIONcaptionsoff
  \newpage
\fi

\begin{IEEEbiography}[{\includegraphics[width=1in,height=1.25in,clip,keepaspectratio]{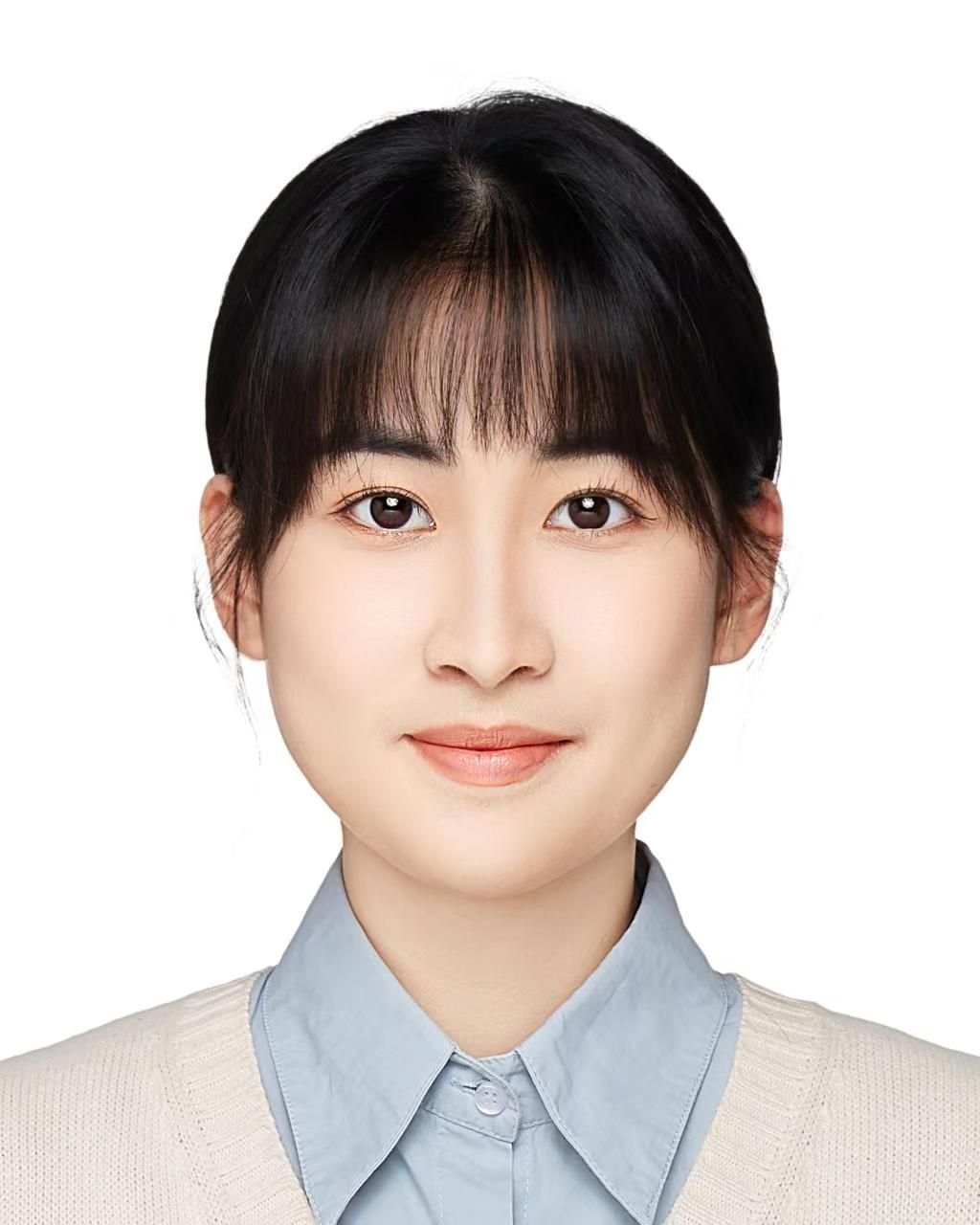}}]
	{Haoyan Lin} received her B.Eng.\ degree in Automation Science and Engineering from South China University of Technology, Guangzhou, China, in 2022. She is currently pursuing the Ph.D.\ degree in the Department of Mechanical and Automation Engineering, The Chinese University of Hong Kong, Hong Kong SAR, China. Her current research interests include Euler-Lagrange systems, output regulation and data-driven control.
	
\end{IEEEbiography}

\begin{IEEEbiography}[{\includegraphics[width=1in,height=1.25in,clip,keepaspectratio]{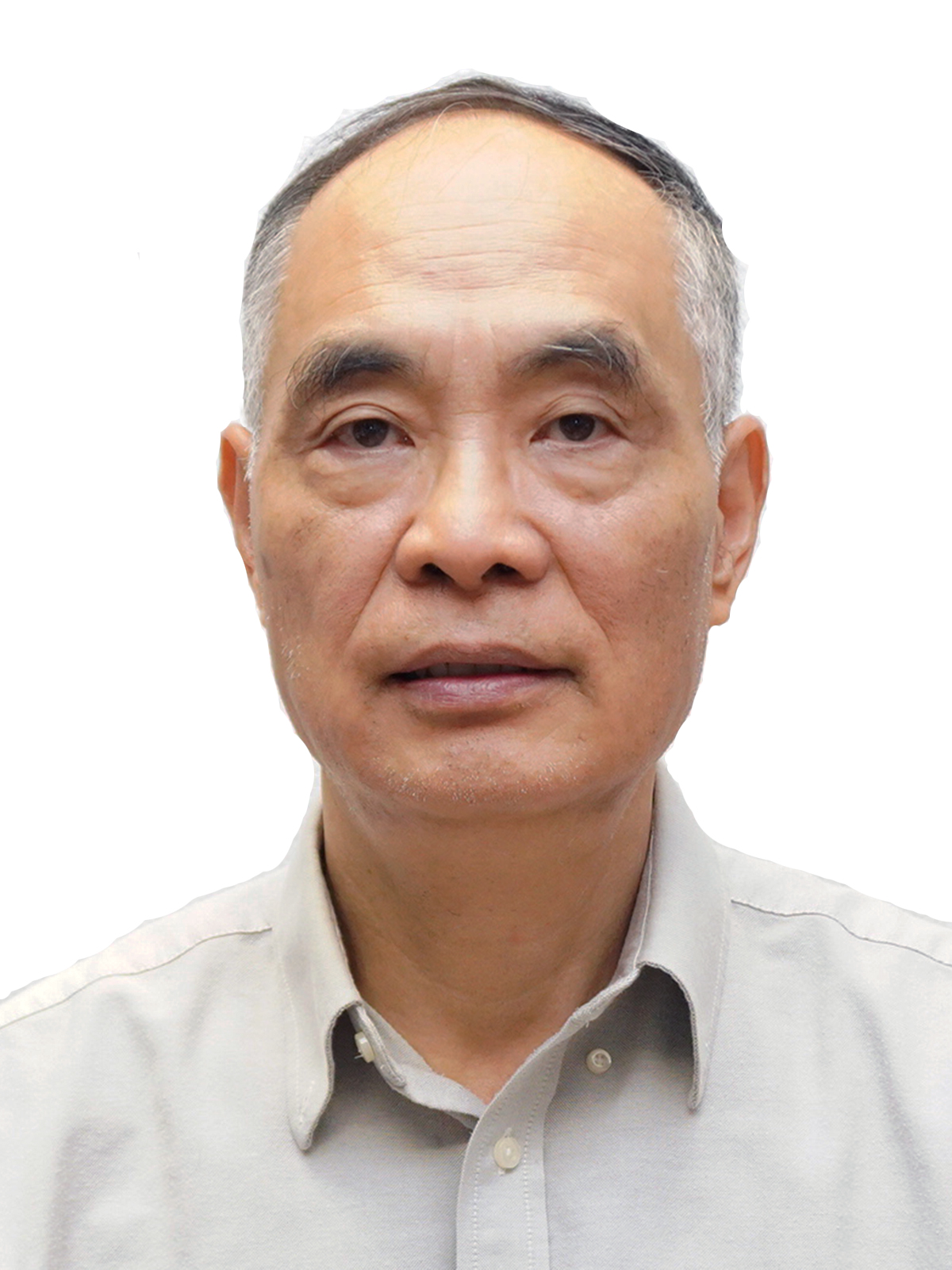}}]{Jie Huang} (Life Fellow, IEEE) received the Diploma
	from Fuzhou University, Fuzhou, China,
	the master’s degree from Nanjing University
	of Science and Technology, Nanjing, China, and
	the Ph.D. degree from Johns Hopkins University,
	Baltimore, MD, USA.
He is currently a Choh-Ming Li Research Professor
of Mechanical and Automation Engineering at The
Chinese University of Hong Kong (CUHK), and Associate Dean (Research) of Faculty of Engineering at CUHK. He is also a Chair Professor with Harbin Institute of Technology. 
His research interests include nonlinear control theory and applications, multi-agent systems, game theory, and flight guidance and control.
Dr. Huang is a fellow of  IFAC, CAA, and HKIE, and an International Member of the Chinese Academy of Engineering.
\end{IEEEbiography}


\begin{thebibliography}{99}
\bibitem{Bertsekas2019}
D. Bertsekas,
\emph{Reinforcement learning and optimal control},
Athena Scientific, 2019.


\bibitem{Bian2016}
T. Bian and Z. -P. Jiang,
``Value iteration and adaptive dynamic programming for data-driven adaptive optimal control design'',
\emph{Automatica},
vol. 71, pp. 348-360, 2016.

%\bibitem{Cai2017}
%H. Cai, F. L. Lewis, G. Hu and J. Huang,
%``The adaptive distributed observer approach to the cooperative output regulation of linear multi-agent systems'',
%\emph{Automatica}, vol. 75, pp. 299-305, 2017.


%\bibitem{Cai2022}
%H. Cai, Y. Su and J. Huang,
%\emph{Cooperative control of multi-agent systems: distributed observer and distributed internal model approaches},
%Springer, 2022.


\bibitem{Chen2022}
C. Chen, L. Xie, K. Xie, F. L. Lewis and S. Xie,
``Adaptive optimal output tracking of continuous-time systems via output-feedback-based reinforcement learning'',
\emph{Automatica},
vol. 146, 110581, 2022.


{\color{black}\bibitem{Davision1976}
E. J. Davison,
``The robust control of a servomechanism problem for linear time-invariant multivariable systems'',
\emph{IEEE Transactions on on Automatic Control},
vol. 21, no. 1, pp. 25-34, 1976.


\bibitem{Francis1976}
B. A. Francis and W. M. Wonham,
``The internal model priciple of control theory'',
\emph{Automatica},
vol. 12,  no. 5,  pp. 457-465, 1976.}

\bibitem{Gao2016}
W. Gao and Z. P. Jiang,
``Adaptive dynamic programming and adaptive optimal output regulation of linear systems'',
\emph{IEEE Transactions on on Automatic Control},
vol. 61, no. 12, pp. 4164-4169, 2016.

\bibitem{Huang2004}
J. Huang,
\emph{Nonlinear output regulation: theory and applications},
SIAM, 2004.


\bibitem{Jiang2012}
Y. Jiang and Z. P. Jiang,
``Computational adaptive optimal control for continuous-time linear systems with completely unknown dynamics'',
\emph{Automatica},
vol. 48, no. 10, pp. 2699-2704, 2012.


\bibitem{kucera1972}
V. Kucera,
``A contribution to matrix quadratic equations'',
\emph{IEEE Transactions on Automatic Control}, vol. 17, no. 3, pp. 344-347, 1972.

{\color{black}\bibitem{Lewis2012}
F. L. Lewis, D. Vrabie, and V. L. Syrmos,
\emph{Optimal Control},
Hoboken, NJ, 547USA: Wiley, 2012.}


\bibitem{LH2024}
L. Lin and J. Huang,
``Data-driven optimal output regulation for continuous-time linear systems via internal model principle'',
\emph{IEEE Transactions on Automatic Control}, 	
vol. 70, no. 6, pp. 4202-4208, 2025.

{\color{black}\bibitem{Lin2024}
L. Lin, H. Lin and J. Huang,
``A new approach to the data-driven output-based LQR problem of continuous-time linear systems'', \emph{IEEE Transactions on Automatic Control}, conditionally accepted. Also available in arXiv preprint arXiv:2509.18819, 2025.}


{\color{black}\bibitem{Liu2018}
Y. Liu and W. Gao,
``Adaptive optimal output regulation of continuous-time linear systems via internal model principle'',
\emph{2018 9th IEEE Annual Ubiquitous Computing, Electronics  \& Mobile Communication Conference (UEMCON)},
pp. 38–43, 2018.


\bibitem{Meyn2022}
S. Meyn,
\emph{Control systems \& reinforcement learning},
Cambridge University Press, 2022.


{\color{black}\bibitem{Modares2016}
H. Modares, F. L. Lewis and Z. P. Jiang,
``Optimal output-feedback control of unknown continuous-time linear systems using off-policy reinforcement learning'',
\emph{IEEE Transactions on Cybernetics}, vol. 46, no. 11, pp. 2401-2410, 2016.}



{\color{black}\bibitem{Rizvi2023}
S. A. A. Rizvi and Z. Lin,
\emph{Output Feedback Reinforcement Learning Control for Linear Systems}, Birkhäuser, 2023.}


\bibitem{Rizvi2019}
S. A. A. Rizvi and Z. Lin,
``Reinforcement learning-based linear quadratic regulation of continuous-time systems using dynamic output feedback'',
\emph{IEEE Transactions on Cybernetics},
vol. 50, no. 11, pp. 4670-4679, 2019.


\bibitem{sutton2018reinforcement}
R. S. Sutton and A. G. Barto,
\emph{Reinforcement learning: An introduction},
MIT press, 2018.


\bibitem{Vrabie2009}
D. Vrabie, O. Pastravanu, M. Abu-Khalaf and F. Lewis,
``Adaptive optimal control for continuous-time linear systems based on policy iteration'',
\emph{Automatica},
vol. 45, no. 2, pp. 477-484, 2009.



\bibitem{Werbos}
P. Werbos,
 ``Neural networks for control and system identification'',
\emph{Proceedings of 1989 IEEE Conference on Decision and Control}, pp. 260-265, 1969.

\bibitem{Wonham}
W. M. Wonham,
``On a matrix Riccati equation of stochastic control'',
\emph{SIAM Journal on Control}, vol. 6, no. 4, pp. 681-697, 1968.


%\bibitem{Su2012}
%Y. Su and J. Huang,
%``Cooperative output regulation of linear multi-agent systems'',
%\emph{IEEE Transactions on on Automatic Control},
%vol. 57, no. 4, pp. 1062-1066, 2012.

%\bibitem{Gao2017}
%W. Gao, Z. P. Jiang, F. Lewis and Y. Wang,
%``Cooperative optimal output regulation of multi-agent systems using adaptive dynamic programming'',
%\emph{American Control Conference (ACC)},
%pp. 2674-2679, 2017.




\bibitem{Xie2022}
K. Xie, X. Yu, and W. Lan,
``Optimal output regulation for unknown continuous-time linear systems by internal model and adaptive dynamic programming'',
\emph{Automatica},
vol. 146, 110564, 2022.}





\bibitem{Xie2023}
K. Xie, Y. Zheng, W. Lan and X. Yu,
``Adaptive optimal output regulation of unknown linear continuous-time systems by dynamic output feedback and value iteration'',
\emph{Control Engineering Practice},
vol. 141, 105675, 2023.





{\color{black}\bibitem{Zhu2014}
L. M. Zhu, H. Modares, G. O. Peen, F. L. Lewis and B. Yue,
``Adaptive suboptimal output-feedback control for linear systems using integral reinforcement learning'',
\emph{IEEE Transactions on Control Systems Technology}, vol. 23, no. 1, pp. 264-273, 2014.}

\end{thebibliography}
\end{document}